\documentclass[12pt]{article}        

 \usepackage[T1]{fontenc}             
 \usepackage[width=16cm, height=22cm]{geometry}   

 \usepackage{mathtools} 
 \usepackage{amssymb}          
\usepackage{hyperref} 
 \usepackage[amsmath,thmmarks,hyperref]{ntheorem} 
 \usepackage{shuffle}

 \usepackage{icomma}                   
 \usepackage{enumerate}          
 \usepackage{color}                 
 \usepackage{setspace}                 
 \usepackage{needspace}               
 \usepackage{url}                    
 \usepackage{mathrsfs}                
 \usepackage{esint}    
 \usepackage{slashed}

 \numberwithin{equation}{section}
 
 \usepackage[numbers,square]{natbib}
 \usepackage{tikz}
\usetikzlibrary{arrows.meta}
\usetikzlibrary{cd}

\theoremstyle{nonumberplain}  
\theoremheaderfont{\itshape}  
\theorembodyfont{\normalfont}  
\theoremseparator{.}  
\theoremsymbol{$\Box$}
\newtheorem{proof}{Proof} 
  
\theoremstyle{plain}  
\theoremheaderfont{\normalfont\bfseries}  
\theorembodyfont{\itshape}  
\theoremsymbol{~}
\theoremseparator{.}  
\newtheorem{proposition}{Proposition}[section]

\newtheorem{lemma}[proposition]{Lemma}  
  
\newtheorem{theorem}[proposition]{Theorem}   

\theorembodyfont{\normalfont}  
 
\newtheorem{remark}[proposition]{Remark}

\newtheorem{definition}[proposition]{Definition}

\theoremstyle{nonumberplain}
\theoremheaderfont{\normalfont\bfseries}  
\theorembodyfont{\itshape}  
\theoremsymbol{~}
\theoremseparator{.}

\newcommand{\R}{\mathbb{R}}
\newcommand{\V}{\mathcal{V}}
\newcommand{\N}{\mathbb{N}}
\newcommand{\Z}{\mathbb{Z}}
\newcommand{\C}{\mathbb{C}}
\newcommand{\dd}{\mathrm{d}}

\newcommand{\ind}{\operatorname{ind}}
\newcommand{\tr}{\mathrm{tr}}
\newcommand{\End}{\mathrm{End}}

\newcommand{\ch}{\mathrm{ch}}
\newcommand{\Ch}{\mathrm{Ch}}

\renewcommand*{\div}{\operatorname{div}}
\newcommand{\cc}{\mathbf{c}}
\newcommand{\cd}{\slashed{\mathbf{c}}}
\DeclareMathOperator{\str}{\mathrm{str}}

\newcommand{\id}{\mathrm{id}}

\renewcommand{\tilde}{\widetilde}

\newcommand{\Dirac}{\mathsf{D}}

\newcommand{\B}{\mathsf{B}}

\newcommand{\T}{\mathbb{T}}

\renewcommand{\L}{\text{\normalfont\sffamily L}}

\newcommand{\sgn}{\mathrm{sgn}}
\newcommand{\Spin}{\mathrm{Spin}}

\newcommand{\Cl}{\mathrm{Cl} }

\DeclareMathOperator{\Str}{\mathrm{Str}}

\newcommand{\defeq}{\stackrel{\mathrm{def}}{=}}

\newcommand{\xx}{\mathsf{x}}
\newcommand{\UU}{\mathsf{u}}
\newcommand{\INT}{\mathcal{I}}
\newcommand{\qq}{\mathsf{q}}

\newcommand{\bfpsi}{\boldsymbol{\psi}}
\newcommand{\bfomega}{\boldsymbol{\omega}}
\newcommand{\bfQ}{\mathbf{Q}}
\newcommand{\bfV}{\mathbf{V}}
\newcommand{\bfN}{\mathbf{N}}

\title{A Rigorous Construction of the Supersymmetric Path Integral Associated to a Compact Spin Manifold}
\author{Florian Hanisch\footnote{Universit\"at Potsdam. fhanisch@math.uni-potsdam.de} ~and Matthias Ludewig\footnote{Universit\"at Regensburg, matthias.ludewig@mathematik.uni-regensburg.de}}

\begin{document}


\maketitle

\begin{abstract}
We give a rigorous construction of the path integral in $\mathcal{N}=1/2$ supersymmetry as an integral map for differential forms on the loop space of a compact spin manifold. It is defined on the space of differential forms which can be represented by extended iterated integrals in the sense of Chen and Getzler-Jones-Petrack. Via the iterated integral map, we compare our path integral to the non-commutative loop space Chern character of G\"uneysu and the second author. Our theory provides a rigorous background to various formal proofs of the Atiyah-Singer index theorem using supersymmetric path integrals, as investigated by Alvarez-Gaum\'e, Atiyah, Bismut and Witten.
\end{abstract}

\section{Introduction}

Recently, B.\ G\"uneysu and the second-named author constructed a non-commutative Chern character for Fredholm modules over differential graded algebras \cite{GueneysuLudewig}. Applying this construction to the dg algebra $\Omega(X)$ of differential forms on a compact spin manifold, this is a functional $\Ch_{\Dirac}$ on the cyclic chain complex of $\Omega(X)$, which is given by a formula closely resembling that of the JLO-cocycle \cite{JLO} for Connes' non-commutative Chern character \cite{ConnesNCDG, ConnesMoscoviciChern}, but with further correction terms coming from the fact that the action of $\Omega(X)$ by Clifford multiplication on the Hilbert space of $L^2$-sections of the spinor bundle is not multiplicative; see formula \eqref{QuantizationMap} below.

On the other hand, in supersymmetric quantum mechanics, one is interested in the path integral of the $\mathcal{N}=1/2$ supersymmetric $\sigma$-model associated to a closed spin manifold $X$. This is an integral over the supermanifold of maps $S^{1|1} \rightarrow X$ which, when the supergeometry is translated into the language of differential forms, can be written as the differential form integral \cite{Atiyah, Lott}
\begin{equation} \label{PathIntegralIntro}
  \INT(\theta) \stackrel{\text{formally}}{=} \int_{\L X} e^{-S-\omega} \wedge \theta, 
\end{equation}
over the ordinary smooth loop space $\L X$, where $S$ is the energy functional and $\omega$ is the canonical two-form on $\L X$ (see \eqref{CanonicalTwoForm} below), exponentiated in the exterior algebra. 

As there is no integration theory for differential forms in infinite dimensions, it is unclear how to give a rigorous meaning to the right hand side of \eqref{PathIntegralIntro}. However, the task of giving a rigorous construction of this {\em supersymmetric path integral} has received a lot of attention ever since Atiyah \cite{Atiyah} and Bismut \cite{Bismut1, BismutLoc, BismutDH} used it formally as a tool to give a ``proof'' of the Atiyah-Singer index theorem. A similar ``proof'' was given by Alvarez-Gaum\'e \cite{AlvarezGaume}, based on ideas of Witten \cite{Witten2, Witten}; the path integral \eqref{PathIntegralIntro} is the transfer of their supergeometric formulation to the language of differential forms. 

\medskip

The purpose of this paper is to clarify the connection between the formal expression \eqref{PathIntegralIntro} for the path integral and the Chern character $\Ch_{\Dirac}$ and to argue that, in fact, the Chern character ``is'' precisely the desired path integral, or rather, more precisely, its pushforward along Chen's \emph{iterated integral map} \cite{Chen1, GJP}.  Using the results of \cite{HanischLudewig2} and the Wiener measure, we give a stochastic interpretation of the formula \eqref{PathIntegralIntro} and compare this path integral with the Chern character. 
Together with \cite{GueneysuLudewig} and \cite{HanischLudewig2}, this paper completes the program (which was to our knowledge initially started by Getzler \cite{GJP0, GetzlerThom, GetzlerOdd, GJP}) to rigorously construct the path integral \eqref{PathIntegralIntro} using the iterated integral map and the cyclic cohomology of $\Omega(X)$.

\paragraph{Domain of the path integral.} Our differential form integral will be defined on the subspace $\Omega_{\mathrm{int}}(\L X) \subset \Omega(\L X)$ of \emph{iterated integrals}. These were first introduced by Chen \cite{Chen1} and further developed by Getzler-Jones-Petrack \cite{GJP} as the image of the \emph{(extended) iterated integral map}
\begin{equation*}
  \rho : \mathsf{B}\bigl(\Omega_{\T}(X)\bigr) \longrightarrow \Omega_{\mathrm{int}}(\L X) \subset \Omega(\L X),
\end{equation*}
 the definition of which we recall in \S\ref{SectionIteratedIntegrals}. Here $\mathsf{B}(\Omega_{\T}(X))$ is the bar complex of the algebra $\Omega_\T(X) = \Omega(X \times\T)^{\T}$, (with $\T = \R/\Z$), which is an enlargement of $\Omega(X)$ needed in order to include the Bismut-Chern character forms (first introduced by Bismut \cite{Bismut1} and Wilson \cite{Wilson} in the odd case) into the theory \cite{CacciatoriGueneysu, GJP}. 
 {Elements $\vartheta \in \Omega_{\T}(X)$ will always be written as $\vartheta = \vartheta^\prime + dt \wedge \vartheta^{\prime\prime}$ with $\vartheta^\prime, \vartheta^{\prime\prime} \in \Omega(X)$.}
It turns out that the Chern character vanishes on the kernel of the iterated integral map 
 and hence can be pushed forward to a functional $\rho_!\Ch_{\Dirac}$ on $\Omega_{\mathrm{int}}(\L X)$.By the properties of the Chern character, it is equivariantly closed (this is the supersymmetry property) and satisfies a localization formula; see \eqref{LocalizationFormula} below.
 
\paragraph{Definition of the path integral.} In order to make sense of the formal expression on the right hand side of \eqref{PathIntegralIntro}, in a first step, one has to make sense of the \emph{top degree part}, or \emph{Berezinian} of the differential form $e^{-\omega} \wedge \theta$ for any iterated integral $\theta$. This is clearly problematic by infinite-dimensionality of the loop space. 

Our starting point here is a formula from our paper \cite{HanischLudewig2}, which  for any spin manifold $X$ and any iterated integral $\theta \in \Omega_{\mathrm{int}}(\L X)$ provides a good interpretation of the top degree of the \emph{composite} form $e^{-\omega} \wedge \theta$. Explicitly, if $\theta = \rho(\vartheta_1, \dots, \vartheta_N)$ with $\vartheta_1, \dots, \vartheta_N \in \Omega_{\T}(X)$, the formula is
\begin{equation} \label{ExplicitFormulaqIntro}
 \int_{\Delta_N}  [\gamma\|_{\tau_1}^0]^\Sigma \prod_{a=1}^N \Bigl(\cd\bigl(\vartheta_a^{\prime\prime}(\gamma_{\tau_a})\bigr) -\cd\bigl(\iota_K \vartheta_a^\prime(\gamma_{\tau_a})\bigr)\Bigr) [\gamma\|_{\tau_{a+1}}^{\tau_a}]^\Sigma \dd \tau,
\end{equation}
where $[\gamma\|_{\bullet}^\bullet]^\Sigma$ denotes parallel transport along the loop $\gamma \in \L X$ in the spinor bundle, $\cd$ denotes (rescaled) Clifford multiplication, $\iota_K$ denotes insertion of the velocity vector field  $K(\gamma) = \dot{\gamma}$ and $\Delta_N = \{0 \leq \tau_1 \leq \dots \leq \tau_N \leq 1\}$ denotes the $N$-dimensional simplex. 

Using the stochastic parallel transport and (iterated) Stratonovich integrals, the expression \eqref{ExplicitFormulaqIntro} has a canonical interpretation as the supertrace of a Clifford-algebra-valued function $Q(\theta)$, which is measurable with respect to the Wiener measure. Since the Wiener measure on the loop space is formally given by $e^{-S} \dd \gamma$, we arrive at the path integral formula
\begin{equation} \label{DefPathIntegralIntro}
  \INT(\theta) ~\defeq~ \int_X \mathbb{E}_x\left[ \left.\str Q(\theta)_{\xx_\bullet} \exp\left(-\frac{1}{8}\int_0^1 \mathrm{scal}(\xx_t) \dd t \right) \right| \xx_1 = x \right] \dd x,
\end{equation}
 where the conditional expectation is taken over a Brownian motion $\xx_\bullet$ with $\xx_0 = \xx_1 = x$. The appearance of scalar curvature term not present in \eqref{PathIntegralIntro} is somewhat of a curiosity, owing to the Lichnerowicz formula; see Remark~\ref{RemarkScalarCurvature}.

\paragraph{Relation to the Chern character.} As mentioned above, the Chern character $\Ch_{\Dirac}$, initially defined on $\mathsf{B}(\Omega_{\T}(X))$, can be pushed forward to a functional $\rho_! \Ch_{\Dirac}$ on $\Omega_{\mathrm{int}}(\L X)$. Ideally, one would like to prove that this pushforward is precisely the path integral \eqref{DefPathIntegralIntro}, in other words the equality
\begin{equation} \label{EqualityICh}
   \rho_!\Ch_{\Dirac} = \INT ~~:~~ \Omega_{\mathrm{int}}(\L X) \longrightarrow \R.
\end{equation}
However, it turns out that this is only true on a certain subset of $\Omega_{\mathrm{int}}(\L X)$; see Thm.~\ref{TheoremEqualityChPI}; the failure for \eqref{EqualityICh} to hold for all entire iterated integrals seems to be due to the bad interaction of certain singularities of the fermionic integral with stochastic integration. 
The ``correct'' rigorous version of the right hand side of \eqref{PathIntegralIntro} is certainly $\rho_! \Ch_{\Dirac}$, {as this satisfies ``supersymmetry'', i.e., is equivariantly closed.}

In any event, our results reveal a remarkable connection between cyclic homology and stochastic analysis \cite{HuMeyer1}. More specifically, the complicated combinatorial structure of the Chern character $\Ch_{\Dirac}$ coming from the failure of $\cc: \Omega(X) \rightarrow \Cl(X)$ to be an algebra homomorphism (see formula \eqref{QuantizationMap} below) has its stochastic counterpart in so-called \emph{Hu-Meyer formulas} that arise when converting iterated Stratonovich integrals into their It\^o version. 

\medskip

\paragraph{Applications.} Incidentally, all ``interesting'' forms seem to lie in the subset of $\Omega_{\int}(\L X)$ where $\INT$ coincides with the Chern character; in particular, this is true for all Bismut-Chern character forms (reviewed in \S\ref{SectionBismutChernCharacters}). Their integrals can be calculated both on the stochastic side and on the chain complex side. Explicitly, we have the formula
\begin{equation*}
  \Ch_{\Dirac} \bigl(\Ch(E)\bigr) = i^{n/2} \ind(\Dirac_E),
\end{equation*}
where $\Ch(E)$ is the Bismut-Chern character of a vector bundle with connection $E$ and $\Dirac_E$ is the Dirac operator twisted with this bundle. The integral of an odd Chern character of Wilson \cite{Wilson} can be expressed in terms of a spectral flow; see Prop.~\ref{PropFormulasBC}.

The localization formula \eqref{LocalizationFormula} for the Chern character $\Ch_{\Dirac}$ translates to a formula of Duistermaat-Heckmann type \cite{DuistermaatHeckmann, BerlineVergne} for the path integral $\INT$. Applying this to the Bismut-Chern characters, one obtains the Atiyah-Singer index formula for the twisted Dirac operator, as envisioned by Atiyah, Bismut and others.

\medskip

A formula similar to \eqref{DefPathIntegralIntro} has previously been given by Lott \cite[\S V]{Lott}; in fact, his paper was one of the main inspirations for our considerations. However, his formula only works for a small subset of our $\Omega_{\mathrm{int}}(\L X)$ and he does not investigate any connections to the cyclic complex.

We would also like to remark that a construction of certain supersymmetric path integrals has been given by Fine and Sawin \cite{FineSawin2, FineSawin1}, using finite-dimensional approximation. From our point of view, they construct the path integral $\INT$ for specific integrands, essentially the Bismut-Chern forms  considered in \S\ref{SectionBismutChernCharacters}.

\medskip

After recalling some basic facts regarding loop space differential forms and spin geometry, we define the bar complex, which is the home for both the iterated integral map and the Chern character $\Ch_{\Dirac}$, which are introduced next. Then in \S\ref{SectionPathIntegral}, we introduce the path integral map $\INT$ and compare it to the Chern character $\Ch_{\Dirac}$. Finally, in \S\ref{SectionBismutChernCharacters}, we recall the definition of the Bismut-Chern characters and calculate their path integrals.

\paragraph{Acknowledgements.} 
We are indebted to J.-M.\ Bismut and B.\ G\"uneysu for helpful discussions.
We thank the Max-Planck-Institute for Gravitational Physics in Potsdam (Albert-Einstein-Institute), the Max-Planck-Institute for Mathematics in Bonn, the Institute for Mathematics at the University of Potsdam and the University of Adelaide for hospitality and financial support.  The second-named author was supported by the Max-Planck-Foundation and the ARC Discovery Project grant FL170100020 under Chief Investigator and Australian Laureate Fellow Mathai Varghese.

\tableofcontents

\section{Preliminaries}

In this preliminary section, we collect some basic facts regarding differential forms on the loop space of a manifold and spin geometry.

\subsection{Loop space differential forms} \label{SectionLoopSpaceForms}

In this section, we give a brief overview of the theory of differential forms on the loop space $\L X = C^\infty(\T, X)$ of a Riemannian manifold $X$, where throughout, we denote $\T = \R/\Z$. This is an infinite-dimensional manifold, modelled on the Fr\'echet space $C^\infty(\T, \R^n)$, where $n = \dim(X)$. Its tangent space $T_\gamma X$ at $\gamma \in \L X$ is canonically identified with $C^\infty(\T, \gamma^*TX)$, the space of sections of the pullback bundle $\gamma^* TX$ over $\T$.

Before we discuss the general definition of differential forms on $\L X$, let us start with some examples. Given a differential form $\vartheta \in \Omega^N(X)$ and $t \in \T$, setting
\begin{equation} \label{TrivialLift}
  \vartheta(t)_\gamma[V_1, \dots, V_N] \stackrel{\text{def}}{=} \vartheta_{\gamma(t)}\bigl[V_1(t), \dots, V_N(t)\bigr]
\end{equation}
for $V_1, \dots, V_N \in T_\gamma \L X$ produces a differential form on $\Omega^N(\L X)$ that we denote by $\vartheta(t)$. Many more examples of differential forms $\L X$ building on this are given below in \S\ref{SectionIteratedIntegrals}, see in particular \eqref{ExtendedIteratedIntegralMap}.

\paragraph{Differential forms on infinite-dimensional manifolds.} 
Let $Y$ be smooth manifold, modelled on a (possibly infinite-dimensional) locally convex space. It turns out that the ``correct'' definition of differential $\ell$-forms on $Y$ is to set 
\begin{equation} \label{DefinitionOfFormsInGeneral}
\Omega^\ell(Y) \defeq C^\infty\bigl(Y, L^\ell_{\mathrm{alt}}(TY, \R)\bigr),
\end{equation}
the space of smooth sections of the vector bundle  $L^\ell_{\mathrm{alt}}(TY, \R)$ over $Y$, the fiber of which at $y \in Y$ is the space of bounded, alternating multilinear functionals on $T_y Y$ \cite[\S33]{KrieglMichor}. Here one uses the notion of {\itshape convenient} smoothness, meaning that smooth curves are mapped to smooth curves, a notion extensively discussed in \cite{KrieglMichor}.
 With this definition, one has a well-defined wedge product, exterior differential, pullback maps and Lie derivatives, just as in finite dimensions (in contrast to several other possible definitions, e.g.\ sections of the exterior power $\Lambda^\ell T^* Y$ of the cotangent bundle, see  \cite[33.21]{KrieglMichor}). 
Setting $Y = \L X$, we let
\begin{equation} \label{LoopSpaceForms}
  \Omega(\L X) \defeq \bigoplus_{\ell=0}^\infty \Omega^\ell(\L X).
\end{equation}
be the dg algebra of differential forms on $\L X$. 

\paragraph{Diffeological differential forms.} To better understand the notion of smoothness present in \eqref{DefinitionOfFormsInGeneral}, we now introduce a second notion of differential forms, which happens to coincide with the previous one for the loop space. This notion of differential forms comes from  viewing an infinite-dimensional manifold $Y$ as a \emph{diffeological space}. Diffeological spaces are given in terms of plots (for details on diffeological spaces, see e.g.\ \cite{ZemmourDiffeology}). If $Y$ already has a manifold structure, a collection of plots is given taking smooth maps $f: S\rightarrow Y$ for $S \subset \R^m$ or, more generally, $S$ any finite-dimensional manifold. A differential form $\theta$ on $Y$ then consists by definition of a collection of differential forms $\theta_f \in \Omega(S)$ for each plot $f: S \rightarrow Y$, subject to the compatibility condition $\theta_{f_1} = g^* \theta_{f_2}$ whenever $f_i: S_i \rightarrow Y$ are two plots such that $f_1 = f_2 \circ g$ for  some smooth map $g: S_1 \rightarrow S_2$. We will call these \emph{diffeological differential forms}. We remark that for $Y = \L X$ and $S$ finite-dimensional, a map $f: S \rightarrow \L X$ is smooth (in the sense that it maps smooth curves in $S$ to smooth curves in $\L X$) if and only if the corresponding map $f^\vee: S \times \T \rightarrow X$ is smooth.

Clearly, any differential form $\theta \in \Omega^\ell(Y)$ gives rise to a diffeological form by setting $\theta_f = f^*\theta$ for a smooth map $f: S \rightarrow Y$; the compatibility follows from naturality of the pullback and one easily checks that the assignment $\theta \mapsto \{\theta_f\}_{f \in \mathrm{Plots}}$ is injective.

\begin{lemma} \label{LemmaDiffeologicalFormsEqualOrdinaryOnes}
For the loop space $\L X$, any diffeological form $\{\theta_f\}_{f \in \mathrm{Plots}}$ is of the form $\theta_f = f^*\theta$ for some $\theta \in \Omega(\L X)$.
\end{lemma}

\begin{proof}
Let $\{\theta_f\}_{f\in\text{Plots}}$ be a diffeological differential $\ell$-form on $\L X$. For $\gamma \in \L X$ and non-zero $V_1, \dots, V_\ell \in T_\gamma Y$, choose a plot $f: \R^\ell \rightarrow Y$ such that $df_0(e_i) = V_i$, $i=1, \dots, \ell$ and set 
\begin{equation*}
\theta_\gamma[V_1, \dots, V_\ell] \defeq \theta_f[e_1, \dots, e_n].
\end{equation*} 
Using the compatibility of the family $\{\theta_f\}$, one shows that this definition is independent of the choice of plot.

So far, for each $\gamma \in \L X$, we have constructed an element $\theta_\gamma \in L^\ell_{\mathrm{alt}}(T_\gamma \L X, \R)$. This gives the desired section of $L^\ell_{\mathrm{alt}}(T \L X, \R)$, but we have to verify its smoothness. In a local parametrization $\kappa : U \rightarrow \L X$, where $U \subset E = C^\infty(\T, \R^n)$ is an open set of the model space, $\kappa^*\theta$ is a map $U \rightarrow L^\ell_{\mathrm{alt}}(E, \R)$; smoothness can now be verified for each such parametrization $\kappa$. To this end, we have to verify that for each curve $c: \R \rightarrow U$, $\kappa^*\theta \circ c$ is a smooth curve in $L^\ell_{\mathrm{alt}}(E, \R)$. 
Let $c$ be such a curve.
Then for any $a \in \R$ and $V_1, \dots, V_\ell \in E$, there exists $\varepsilon >0$ such that $c(r) + s_1 V_1 + \dots + s_\ell V_\ell \in U$ for all $s_i \in [-\varepsilon, \varepsilon]$ and $r \in [-a, a]$. We then define a plot $f: (-\varepsilon, \varepsilon)^\ell \times (-a, a) \longrightarrow \L X$ by
\begin{equation*}
  f(s_1, \dots, s_\ell, r) = \kappa \bigl( s_1 V_1 + \dots + s_\ell V_\ell + c(r)\bigr).
\end{equation*}
Then 
\begin{equation*}
  (\kappa^* \theta \circ c)(r)[V_1, \dots, V_\ell] = \theta_{f}|_{0, \dots, 0, r}\bigl[ e_1, \dots, e_\ell\bigr],
\end{equation*}
which is a smooth function of $r$, since $\theta_f$ is a smooth differential form on $(-\varepsilon, \varepsilon)^\ell \times (-a, a)$. The smoothness of $(\kappa^*\theta \circ c)(r)$ as a function of $r \in (-a, a)$ now follows from  Thm.~5.18 in \cite{KrieglMichor}, which states that smoothness can be tested pointwise in our case.
\end{proof}

The space $\Omega(\L X)$ of differential forms   on the loop space has a natural topology coming from considering $\L X$ as a diffeological space. 
It is the initial topology induced by the pullback maps $f^* : \Omega(\L X) \rightarrow \Omega(S)$ for smooth maps $f:S \rightarrow \L X$ (such a map is smooth if and only if the corresponding map $f^\vee: S \times \T \rightarrow X$ is smooth). 
This is a (degreewise) complete locally convex topology on $\Omega(\L X)$; a family of continuous seminorms inducing this topology is the family $\nu_f(\theta) = \nu(f^*\theta)$, where $f$ ranges over all smooth maps from smooth manifolds $S$ into $\L X$ and $\nu$ is a continuous seminorm on $\Omega(S)$. 
To see completeness, let $\theta_i \in \Omega(\L X)$, $i \in \mathscr{I}$ be a Cauchy net. Equivalently, this means that all the nets $f^*\theta_i$ are Cauchy, for all $f: S \rightarrow \L X$, hence by completeness of $\Omega(S)$, we obtain limits $\theta_f \in \Omega(S)$. 
Whenever $f_k: S_k \rightarrow \L X$ ($k=1, 2$) and $g: S_1 \rightarrow S_2$ are such that $f_2 \circ g = f_1$, we have $g^* f_2^* \theta_i = f_1^*\theta_i$ and hence in the limit, we obtain $g^*\theta_{f_2} = \theta_{f_1}$.
We conclude that the collection $\{\theta_f\}_{f \in \mathrm{Plots}}$ defines a diffeological differential form $\theta$ on $\L X$, which is the same as an ordinary form by Lemma~\ref{LemmaDiffeologicalFormsEqualOrdinaryOnes}.

\paragraph{The circle action.} So far, we ignored an important feature of the free loop space, which is its circle action given by rotation of loops,
\begin{equation*}
~~~~~~~~~~~~~~~~~~~~~~~~~~~~~~~ (\tau \cdot \gamma)(t) = \gamma(t+\tau), \qquad t, \tau \in \T = \R/\Z.
\end{equation*}
The generating vector field for the circle action is the vector field
\begin{equation} \label{KillingField}
K(\gamma) \defeq \dot{\gamma} \in T_\gamma \L X,
\end{equation}
where the notation references to the fact that it is a Killing vector field with respect to the canonical $L^2$-metric on $\L X$, generating a one-parameter family of isometries. These data induce the {\em equivariant} differential
\begin{equation}\label{EquivariantDifferential}
d_K \defeq d - \iota_K,
\end{equation}
where $\iota_K$ denotes insertion of the vector field $K$.  By Cartan's formula, $d_K$ squares to the Lie derivative $\mathscr{L}_K$ with respect to $K$, hence $d_K$ is a differential on the complex $\Omega(\L X)^\T$ of $\T$-invariant differential forms, which is $\Z_2$-graded by reducing the usual grading mod $2$ (notice that $d_K$ is inhomogeneous, but does preserve the $\Z_2$-grading).

Throughout, $\Omega(\L X)$ denotes the algebra of differential forms that are a direct \emph{sum} of its homogeneous components; see \eqref{LoopSpaceForms}.
It is well known, however, that in equivariant cohomology of the loop space, it is important to allow differential forms that are an infinite sum of its homogeneous components, in other words, elements of the direct {\em product} of the homogeneous components \cite{JonesPetrack}. Denote the direct product of the even respectively odd forms by
\begin{equation} \label{OmegaHat}
  \widehat{\Omega}^{+}(\L X) \stackrel{\mathrm{def}}{=} \prod_{\ell=0}^\infty \Omega^{2\ell}(\L X), \qquad \widehat{\Omega}^{-}(\L X) \stackrel{\mathrm{def}}{=} \prod_{\ell=0}^\infty \Omega^{2\ell+1}(\L X).
\end{equation}
The $\T$-invariant subspace $\widehat{\Omega}(\L X)^{\T} \subset \widehat{\Omega}(\L X)$ is a $\Z_2$-graded complex with respect to the equivariant differential $d_K$ \eqref{EquivariantDifferential}. This complex will be of importance in \S\ref{SectionBismutChernCharacters}.

The action functional and the canonical two-form on $\L X$ appearing in the path integral formula \eqref{PathIntegralIntro} are explicitly given by
\begin{equation} \label{CanonicalTwoForm}
 S(\gamma) = \frac{1}{2} \int_{\T} |\dot{\gamma}(t)|^2 \dd t, \qquad \omega[V, W] \defeq \int_{\T} \bigl\langle V(t), \nabla_{\dot{\gamma}} W(t)\bigr\rangle\, \dd t, 
\end{equation}
where $V, W \in T_\gamma \L X = C^\infty(\T, \gamma^*TX)$ are smooth vector fields along $\gamma$. The inhomogeneous differential form $S+\omega$ appearing in \eqref{PathIntegralIntro} is closed with respect to the equivariant differential (see \cite[\S V]{Lott} or \cite{Atiyah, Bismut1}).

%

\subsection{Spin Geometry} \label{SectionSpinPreliminaries}

In this section, we give a quick account of some notions of spin geometry needed in this paper. 
Let $X$ be a Riemannian manifold. 

\paragraph{The real spinor bundle.} A \emph{spin structure} is a lift of the structure group form $O(n)$ to $\Spin_n$. In particular, this is a $\Spin_n$-principal bundle $P$. The associated {\em (real) spinor bundle} is
\begin{equation} \label{AssociatedSpinorBundle}
  \Sigma \defeq P \times_{\Spin_n} \Cl_n,
\end{equation}
Here (as usual) we identify $[p\cdot g, a] = [p, g \cdot a]$, $g \in \Spin_n$, inside $P \times \Cl_n$, where the action of $\Spin_n$ on $\Cl_n$ is by left multiplication, after realizing $\Spin_n$ inside the even part of the Clifford algebra.
$\Sigma$ is naturally a bundle of graded $\Cl(TX)$-$\Cl_n$-bimodules on $Y$ (where the grading comes from the even/odd grading of the Clifford algebra). Here $\Cl(TX)$ is the bundle of Clifford algebras over $X$. $\Sigma$ has a canonical connection induced by the Levi-Civita connection on $X$; the right $\Cl_n$-action is parallel with respect to this connection.

Clifford (left-)multiplication by vector fields and its extension to differential forms is denoted by $\cc$. It commutes with the right action of $\Cl_n$. Explicitly, on $\ell$-forms, it is determined by the formula
\begin{equation*}
\cc(e_{i_1} \wedge \cdots \wedge e_{i_\ell}) = \frac{1}{\ell!} \sum_{\sigma \in S_\ell} \sgn(\sigma) \cc(e_{i_{\sigma_1}}) \cdots\cc(e_{i_{\sigma_\ell}}).
\end{equation*}
We also define a rescaled version of Clifford multiplication by
\begin{equation} \label{RescaledCliffordMultiplication}
  \cd(\vartheta) = 2^{-|\vartheta|/2} \cc(\vartheta).
\end{equation}
Here $|\vartheta|$ denotes the degree of $\vartheta$. Throughout, in formulas such as the above, we assume that $\vartheta$ is homogeneous and extend by linearity.
The {\em Dirac operator} $\Dirac$ and its rescaled version $\slashed{\Dirac}$, acting on sections of $\Sigma$, are given by
\begin{equation} \label{FormulaDiracOperator}
\Dirac\psi = \sum_{j=1}^n \cc(e_j) \nabla_{e_j}\psi \qquad \text{and} \qquad \slashed{\Dirac}\psi = \sum_{j=1}^n \cd(e_j) \nabla_{e_j}\psi = 2^{-1/2} \Dirac\psi
\end{equation}
with respect to a local orthonormal basis $e_1, \dots, e_n$ of $TX$. Since the right $\Cl_n$-action is parallel and commutes with Clifford multiplication, the Dirac operator also commutes with the right $\Cl_n$-action.

\paragraph{Supertraces.}
We denote by $\End_{\Cl_n}(\Sigma)$ the space of endomorphisms of the real spinor bundle $\Sigma$ that commute with the right action of $\Cl_n$. Such an endomorphism is always given by left multiplication by an element $a$ of the Clifford algebra $\Cl(TX)$, and we define its {\em supertrace} by the formula
\begin{equation} \label{DefinitionSupertrace}
  \str(a) = 2^{n/2} \langle a, \cc(e_1) \cdots \cc(e_n)\rangle = \langle \cd^{-1}(a),  e_1 \wedge \cdots \wedge e_n\rangle
\end{equation}
for  an \emph{oriented} orthonormal basis $e_1, \dots, e_n$ of $TX$. This is indeed a supertrace, in the sense that it satisfies the graded cyclic permutation property
\begin{equation} \label{CyclicPermutation}
   \str(ab) = (-1)^{|a||b|}\str(ba)
\end{equation}
on homogeneous elements $a, b \in \Cl(TX)$. The supertrace is even, respectively odd (meaning that it vanishes on odd respectively even elements), depending on whether the dimension $n$ of $X$ is even or odd.

Let $A$ be a bounded operator on $L^2(X, \Sigma)$ with smooth integral kernel $a(x, y)$ that commutes with the right action of $\Cl_n$, for example the heat kernel $A = e^{-t \Dirac^2}$ of the Dirac operator. Such an operator has a supertrace, defined by integrating the pointwise supertrace \eqref{DefinitionSupertrace} over the diagonal,
\begin{equation} \label{OperatorSupertrace}
\Str(A) = \int_{X} \str a(x, x) \dd x.
\end{equation}
This supertrace satisfies a cyclic permutation property analogous to \eqref{CyclicPermutation}.

\paragraph{Comparison to the complex spinor bundle.} 
We can also form the {\em complex} spinor bundle
\begin{equation*}
  \Sigma_\C \defeq P \times_{\Spin_n} S,
\end{equation*}
where $S$ is the complex spinor space. If $n$ is even, this is a bundle of irreducible graded modules for $\Cl(TX)$, and hence this bundle satisfies $\End(\Sigma_\C) = \Cl(TX) \otimes \C$. We can therefore also take the operator supertrace $\str_\C$ of elements of $\Cl(TX)$, which is related to  \eqref{DefinitionSupertrace} by the formula
\begin{equation} \label{TraceComparisonEven}
  \str_\C(a) = (-i)^{n/2}\str(a)
\end{equation}
for $a \in \End(\Sigma_\C)$. 
In contrast, if $n$ is odd, then $S$ is not graded, and we have the formula
\begin{equation} \label{TraceComparisonOdd}
  \tr_\C(a) = i^{\frac{n+1}{2}} 2^{n-1/2}\str(a) + 2^m \langle a, \mathbf{1}\rangle
\end{equation}
for the endomorphism trace of $a \in \Cl(\V) \subset \End(\Sigma_\C)$.


\section{Iterated integrals and the Chern character}

In this section, we review the iterated integral map and the definition of the Chern character of \cite{GueneysuLudewig}. Both of these naturally live on the bar complex, which we introduce first.

\subsection{The bar complex} \label{SectionBarComplex}

The {\em bar complex} associated to a differential graded algebra $\Omega$ is the graded vector space
\begin{equation} \label{BarComplex}
  \mathsf{B}\bigl(\Omega\bigr) = \bigoplus_{N=0}^\infty \Omega[1]^{\otimes N}.
\end{equation}
Here $\Omega[1]$ equals $\Omega$ as a vector space, 
but with degrees shifted by one; explicitly, $\vartheta \in \Omega^{k+1}$ iff $\vartheta \in \Omega[1]^{k}$.
{$\Omega[1]^{\otimes N}$ then carries the tensor product grading.}
The elements of $\B(\Omega)$ are called bar chains. It is customary to denote such elements by $(\vartheta_1, \dots, \vartheta_N)$, suppressing the tensor product sign in notation for convenience. 
 $\B(\Omega)$ has two differentials, a differential coming from the differential of $\Omega$ and the {\em bar differential}; they are given by
\begin{equation*}
\begin{aligned}
b_0(\vartheta_1, \dots, \vartheta_N) &= \sum_{k=1}^N (-1)^{n_{k-1}}({\vartheta}_1, \dots, {\vartheta}_{k-1}, d \vartheta_k, \dots, \vartheta_N)\\
b_1(\vartheta_1, \dots, \vartheta_N) &= -\sum_{k=1}^{N-1} (-1)^{n_k}({\vartheta}_1, \dots, {\vartheta}_{k-1}, {\vartheta}_{k}\vartheta_{k+1}, \vartheta_{k+2}, \dots, \vartheta_N),
\end{aligned}
\end{equation*}
where $n_k = |\vartheta_1|+\dots + |\vartheta_k|-k$. 
 The above differentials satisfy $b_0b_1 + b_1 b_0 = 0$, hence turn $\B(\Omega)$ into a bicomplex with total differential $b := b_0 + b_1$.  The differentials moreover descend to the subspace
\begin{equation} \label{CyclicChains}
\mathsf{B}^{\natural}(\Omega) = \mathrm{span}~\Bigl\{ \sum_{k=0}^N (-1)^{n_k(n_N-n_k)} (\vartheta_{k+1}, \dots, \vartheta_N, \vartheta_1, \dots, \vartheta_k) \Bigr\} 
\end{equation} 
of \emph{cyclic chains}, turning it into a subcomplex.

Of course, we can take $\Omega = \Omega(X)$ in the above, the dg algebra of differential forms on a manifold $X$. For our purposes, we more generally consider
\begin{equation*}
  \Omega_{\T}(X) \stackrel{\text{def}}{=} \Omega(X \times \T)^{\T},
\end{equation*}
the space of differential forms on $X \times \T$ which are constant in the $\T$-direction. Elements $\vartheta \in \Omega_\T(X)$ will always be written as
\begin{equation} \label{DecompositionPrime}
\vartheta = \vartheta^\prime + dt \wedge \vartheta^{\prime\prime}, \qquad \text{with} \quad \vartheta^\prime, \vartheta^{\prime\prime} \in \Omega(X).
\end{equation}
The differential of $\Omega_\T(X)$ is $d_\T = d- \iota_{\partial_t}$, where $\iota_{\partial_t}$ denotes insertion of the canonical vector field $\partial_t$ on the $\T$ factor and $d$ denotes the de-Rham differential on $X \times \T$. In terms of \eqref{DecompositionPrime}, we have
\begin{equation*}
  d_\T\vartheta = d_\T(\vartheta^\prime + dt \wedge \vartheta^{\prime\prime}) = d \vartheta^\prime - dt \wedge d\vartheta^{\prime\prime} -  \vartheta^{\prime\prime},
\end{equation*} 
where on the right hand side, $d$ denotes the de-Rham differential on $X$. Observe that this differential is not homogeneous, so that $\Omega_\T(X)$ is only a $\Z_2$-graded complex; however, this does not cause any complications for the theory (in \cite{GueneysuLudewig}, the approach was taken to set $\sigma = dt$ and to declare it to be of degree $-1$).

\paragraph{Entire chains.} In order to deal with the Bismut-Chern characters in \S\ref{SectionBismutChernCharacters}, we need the larger complex of \emph{entire chains} $\mathsf{B}_\epsilon(\Omega_{\T}(M))$, which allows certain infinite sums. It is defined as the completion of $\mathsf{B}(\Omega_{\T}(M))$ with respect to the seminorms
\begin{equation} \label{EntireSeminorms}
  \epsilon_k(c) := \sum_{N=0}^\infty \frac{\pi_{k}^{N}(c_N)}{\lfloor N/2\rfloor!} \qquad \text{for} \qquad c = \sum_{N=0}^\infty c_N,
\end{equation}
where $\pi_{k}^{N}$ denotes the ($N$-fold) projective tensor product norm on $\Omega_{\T}(X)[1]^{\otimes N}$, induced by the $C^k$-norm on $\Omega_{\T}(X)$.
The differential $b$ extends to entire chains, making $\mathsf{B}_\epsilon(\Omega_\T(X))$ a chain complex. In total, we have the following hierarchy of chain complexes:
\begin{equation*}
\begin{tikzcd}[column sep=0cm, row sep = 0cm]
 \text{\footnotesize{(entire cyclic chains)}}~~~ & \mathsf{B}_\epsilon^\natural(\Omega_{\T}(X)) & \subset & \mathsf{B}_\epsilon(\Omega_{\T}(X)) &  ~~~ \text{\footnotesize{(entire chains)}}\\
 & \reflectbox{\rotatebox[origin=c]{90}{$\subset$}} & & \reflectbox{\rotatebox[origin=c]{90}{$\subset$}} &\\
 \text{\footnotesize{(cyclic chains)}} & \mathsf{B}^\natural(\Omega_{\T}(X)) & \subset & \mathsf{B}(\Omega_{\T}(X))& ~~\text{\footnotesize{(bar chains)}}
\end{tikzcd}
\end{equation*}
The Bismut-Chern characters $\mathrm{Ch}(q)$ defined in \S\ref{SectionBismutChernCharacters} live in $\mathsf{B}_\epsilon^\natural (\Omega_\T(X))$ but not in $\mathsf{B}(\Omega_{\T}(X))$. Dually, the Chern character $\mathrm{Ch}_{\Dirac}$ defined in \S\ref{SectionChernCharacter} is a linear functional \emph{a priori} defined on $\mathsf{B}(\Omega_{\T}(X))$, which then turns out to satisfy the necessary estimates to extend to the space of entire chains. In particular, $\mathrm{Ch}_{\Dirac}$ can be evaluated on $\mathrm{Ch}(q)$.

\subsection{Iterated Integrals} \label{SectionIteratedIntegrals}

Given a manifold $Y$, Chen's {\em iterated integral map} $\rho_0: \Omega(Y)^{\otimes N} \rightarrow \Omega(\L Y)$  constructs differential forms on the loop space from differential forms on $Y$ \cite{Chen1}. It can be defined as the composition 
\begin{equation} \label{IteratedIntegralMapComposition}
\begin{tikzcd}
\rho_0 ~ : ~ \Omega(Y)^{\otimes N} \subseteq \Omega\bigl(Y^N\bigr) \ar[r, "\mathrm{ev}^*"]  & \Omega\bigl(\L Y \times \Delta_N\bigr) \ar[r, "\int_{\Delta_N}"] & \Omega\bigl(\L Y\bigr).
 \end{tikzcd}
\end{equation}
where $\mathrm{ev}(\gamma, \tau) = (\gamma(\tau_1), \dots, \gamma(\tau_N))$ and $\int_{\Delta_N}$ denotes integration over the fiber $\Delta_N$. 
 Note that by letting $N$ vary, $\rho_0$ is naturally viewed as a map defined on the bar complex $\B(\Omega(Y))$ defined in \eqref{BarComplex}; since the fiber integration in \eqref{IteratedIntegralMapComposition} has degree $-N$, the grading shift in the definition of $\B(\Omega(Y))$ ensures that $\rho_0$ is grading preserving if $\Omega(\L X)$ carries its usual grading by form degree.

\paragraph{Extended iterated integrals.} For our purposes, we need an extension of this, introduced by Getzler, Jones and Petrack \cite{GJP0, GJP}. 
The {\em extended iterated integral map} is the map
\begin{equation} \label{ExtendedIteratedIntegralMapDom}
  \rho = \alpha^* \circ \rho_0 : \B(\Omega_\T(X)) \longrightarrow \Omega(\L X)
\end{equation}
where $\rho_0$ is the (non-extended) iterated integral map \eqref{IteratedIntegralMapComposition} for the manifold $Y = X \times \T$ and $\alpha^*$ denotes pullback along the map $\alpha: \L X \rightarrow \L(X \times \T)$ given by $\alpha(\gamma)(t) = (\gamma(t), -t)$. Explicitly, $\rho$ is given by the formula
\begin{equation} \label{ExtendedIteratedIntegralMap}
 \rho(\vartheta_1, \dots, \vartheta_N) 
 = \int_{\Delta_N} \bigl(\iota_K \vartheta_1^\prime(\tau_1) - \vartheta_1^{\prime\prime}(\tau_1)\bigr) \wedge \cdots \wedge \bigl(\iota_K \vartheta_N^\prime(\tau_N) - \vartheta_N^{\prime\prime}(\tau_N)\bigr) \dd \tau
 \end{equation}
for $\vartheta_1, \dots, \vartheta_N \in \Omega_\T( X)$, where we recall that $K(\gamma) = \dot{\gamma}$ is the canonical velocity vector field, and the forms $\vartheta_a^{\prime}(\tau_a)$, $\vartheta_a^{\prime \prime}(\tau_a)$ are given by \eqref{TrivialLift}. Again, due to the grading shift in the definition of $\B(\Omega_\T(X))$, $\rho$ is degree-preserving. A tedious calculation shows that $\rho$ is a chain map when restricted on entire cylic chains, explicitly,
\begin{equation} \label{IteratedIntegralChainMap}
\rho\bigl(b(c)\bigr)  + d_K\rho(c) = 0
\end{equation}
for all $c \in \mathsf{B}_\epsilon^\natural(\Omega_\T(X))$, where $d_K$ is the equivariant differential \eqref{EquivariantDifferential}. 

\begin{lemma} \label{LemmaIteratedIntegralMapContinuous}
The iterated integral map is degree-wise continuous when $\B(\Omega_\T(X))$ carries the locally convex topology induced by the seminorms \eqref{EntireSeminorms} and $\Omega(\L X)$ the diffeological topology introduced in \S\ref{SectionLoopSpaceForms}.
\end{lemma}

Here $\Omega(\L X)$ carries the diffeological topology introduced in \eqref{SectionLoopSpaceForms}.

\begin{proof}
Let $f: S \rightarrow \L X$ be a smooth map, and let $f^\vee: S \times \T \rightarrow X$ be the corresponding map defined by $f^\vee(s, \tau) = f(s)(\tau)$. Let $\vartheta \in \Omega^\ell(X)$. Then for $s \in S$ and $v_1, \dots, v_\ell \in T_s S$, we have
\begin{equation*}
\begin{aligned}
  f^*\vartheta(\tau)[v_1, \dots, v_\ell] &= \vartheta_{f^\vee(s, \tau)} \bigl[d f^\vee(s, \tau)v_1, \dots, d f^\vee(s, \tau)v_\ell\bigr],\\
  f^*\iota_K\vartheta(\tau)[v_1, \dots, v_{\ell-1}] &= \vartheta_{f^\vee(s, \tau)} \bigl[\dot{f}^\vee(s, \tau), d f^\vee(s, \tau)v_1, \dots, d f^\vee(s, \tau)v_{\ell-1}\bigr],
\end{aligned}
\end{equation*}
where $d f^\vee$ denotes the differential in the $S$ entry and $\dot{f}^\vee$ the $\tau$-derivative. For a compact set $K\subset S$ and $m \in \N$, we then have 
\begin{equation*}
  \|f^*\vartheta(\tau)\|_{m, K} \leq \|df^\vee\|_{m, K}^\ell \|\vartheta\|_{m, K}, \qquad \|f^*\iota_K \vartheta(\tau)\|_{m, K} \leq \|\dot{f}^\vee\|_{m, K}\|df^\vee\|_{m, K}^{\ell-1} \|\vartheta\|_{m, K},
\end{equation*}
where $\| \cdot \|_{m, K}$ denote the $C^m$-norms on $K$.
In view of \eqref{ExtendedIteratedIntegralMap}, we therefore get
\begin{equation*}
\begin{aligned}
\bigl\|f^* \rho(\vartheta_1, \dots, \vartheta_N)\bigr\|_{K, m} &\leq \frac{1}{N!} \|df\|_{K, m}^{|\vartheta_1| + \dots + |\vartheta_N|} \|\vartheta_1\|_{K, m} \cdots \|\vartheta_N\|_{K, m}\\
&=  \frac{1}{N!} \|df\|_{K, m}^{|\vartheta_1| + \dots + |\vartheta_N|} \nu\bigl((\vartheta_1, \dots, \vartheta_N)\bigr),
\end{aligned}
\end{equation*}
where we denoted $\nu = \|\cdot \|_{K, m}$.
Notice that $|\vartheta_1|+\dots + |\vartheta_N| = \mathrm{deg}(\vartheta_1, \dots, \vartheta_N) + N$, when the degree is taken in $\B(\Omega_\T(X))$. Therefore, if $c = \sum_{N=0}^\infty c_N$ is a chain of total degree $\ell$ with $c_N \in \Omega_\T(X)^{\otimes N}$, then
\begin{equation*}
  \|f^* \rho(c)\|_{K, m} \leq \sum_{N=0}^\infty \frac{1}{N!} \|df\|_{K, m}^{\ell + N} \,\nu(c_N) \leq \sup_{N \in \N_0} \left\{\frac{\|df\|_{K, m}^{N+\ell} \lfloor N/2\rfloor!}{N!}\right\} \nu_\epsilon(c),
\end{equation*}
where one observes that the supremum is finite. Since the topology on $\Omega(X)$ is generated by norms of this type, this proves the claim.
\end{proof}

By the above lemma, any fixed degree $\ell \in \N_0$, we may continuously extend the iterated integral map from $\mathsf{B}^\ell(\Omega_{\T}(X))$ to $\mathsf{B}^\ell_\varepsilon(\Omega_\T(X))$. We can therefore make the following definition.

\begin{definition}[Iterated integrals] 
The space of degree $\ell$ \emph{iterated integrals} $\Omega^\ell_{\mathrm{int}}(\L X) \subset \Omega^\ell(\L X)$ is the image of $\mathsf{B}^\ell_\varepsilon(\Omega_\T(X))$ under the iterated integral map.
\end{definition}

\begin{remark}
Observe that since the degree $\ell$ component of $\mathsf{B}(\Omega_{\T}(X))$ is given by
\begin{equation*}
\mathsf{B}^\ell\bigl(\Omega_{\T}(X)\bigr) = \bigoplus_{N=0}^\infty \bigoplus_{k_1 + \dots + k_N = \ell+N} \Omega^{k_1+1}_{\T}(X) \otimes \cdots \otimes \Omega^{k_N+1}_{\T}(X),
\end{equation*}
for each $N \in \N$ and any $\ell \leq N$, there are elements of degree $\ell$ in $\Omega_\T(X)[1]^{\otimes N}$: e.g., for any $\vartheta \in \Omega^1_{\T}(X)$, the chain $(\vartheta, \dots, \vartheta) \in \Omega_{\T}(X)^{\otimes N}$ has degree zero. Hence if we call elements of the latter space \emph{$N$-ary chains}, the completion $\mathsf{B}_\varepsilon(\Omega_\T(X))$ will contain infinite sums of unbounded arity even in each fixed degree $\ell$. In fact, this happens for the Bismut-Chern characters; see Remark~\ref{InfiniteBismutChern}.
\end{remark}

\subsection{The Chern character} \label{SectionChernCharacter}

In this section, we recall the definition of the \emph{loop space Chern character} 
\begin{equation*}
  \Ch_{\Dirac}: \B\bigl(\Omega_\T(X)\bigr) \longrightarrow \R,
\end{equation*}
associated to a compact spin manifold $X$, introduced in \cite{GueneysuLudewig}. The terminology stems from the fact that $\Ch_{\Dirac}$ is the Chern character of a certain Fredholm module (in the sense of non-commutative geometry) on the dg algebra $\Omega(X)$ (respectively $\Omega_\T(X)$), essentially given by the Dirac operator $\Dirac$; for details, we refer to \cite{GueneysuLudewig}. In \S\ref{SectionComparison}, we relate this to the path integral defined above.

\begin{remark}
The construction presented here differs slightly from that in \cite{GueneysuLudewig}. For simplicity, we use the bar complex instead of the cyclic complex; however, these are related by a map (called $\alpha$ in \cite[(5.2)]{GueneysuLudewig}) such that the Chern character on the cyclic complex is just the pullback of that on the bar complex. Also, we use the real spinor bundle instead of the complex one, resulting in our Chern character being real. Finally, our Hilbert space comes with the structure of a right $\Cl_n$-module, and all operators in consideration commute with this action. Thus, we can use the operator supertrace \eqref{OperatorSupertrace} induced by \eqref{DefinitionSupertrace}, which allows to also handle the case that $\dim(X)$ is odd.
\end{remark}

\paragraph{Cochains.}  
Let $\Omega$ be a $\Z_2$-graded algebra. Given a $\Z_2$-graded algebra\footnote{In the case that $\mathcal{L} = \C$, we endow $\C$ with the trivial grading rendering it purely even and just speak of \emph{bar cochains}.} $\mathcal{L}$, an \emph{$\mathcal{L}$-valued bar cochain} over $\Omega$ is a linear map $\ell: \B(\Omega) \rightarrow \mathcal{L}$. Such a cochain can be viewed as a sequence of multilinear maps 
\begin{equation*}
\ell : \underbrace{\Omega \times \cdots \times \Omega}_N \rightarrow \mathcal{L},
\end{equation*}
again denoted by the same letter. In particular, for $N=0$, this is just an element of $\mathcal{L}$, which we denote by $\ell(\emptyset)$, by abuse of notation.
We say that $\ell$ is \emph{even} if it preserves parity and \emph{odd} if it reverses parity.
The standard coalgebra structure on the tensor algebra $\B(\Omega)$ induces a product on $\mathcal{L}$-valued bar cochains, given by
\begin{equation} \label{ProductOfCochains}
 (\ell_1 \ell_2)(\theta_1, \dots, \theta_N) = \sum_{k=0}^N (-1)^{n_k |\ell_2|} \ell_1(\theta_1, \dots, \theta_k) \ell_2(\theta_{k+1}, \dots, \theta_N),
\end{equation}
where $n_k = |\theta_1|+\dots + |\theta_k|-k$.
This product is compatible with the codifferential $\beta$ defined by
\begin{equation*}
  (\beta\ell)(\theta_1, \dots, \theta_N) = - (-1)^{|\ell|} \ell(b(\theta_1, \dots, \theta_N),
\end{equation*}
in the sense that $\beta(\ell_1 \ell_2) = \beta(\ell_1) \ell_2 + (-1)^{|\ell_1|} \ell_1 \beta(\ell_2)$ for all homogeneous cochains $\ell_1, \ell_2$. In other words, $\beta$ is a derivation on the cochain algebra.

\paragraph{The Chern character.}
If $X$ is a compact spin manifold with spinor bundle $\Sigma$, there is an associated cochain $\omega$ over $\Omega_{\T}(X)$ with values in operators on $L^2(X, \Sigma)$, given by
\begin{equation*}
  \omega(\emptyset) = -\slashed{\Dirac}, \qquad \omega(\vartheta) = \cd(\vartheta^\prime), \qquad \text{and} \qquad \omega(\vartheta_1, \dots, \vartheta_k) = 0, \quad k\geq 2,
\end{equation*}
with $\slashed{\Dirac}$ and $\cd$ being the (rescaled) Dirac operator and Clifford multiplication, see \eqref{RescaledCliffordMultiplication}, \eqref{FormulaDiracOperator}. Here as usual, we write $\vartheta \in \Omega_{\T}(X)$ as $\vartheta = \vartheta^\prime + dt \wedge \vartheta^{\prime\prime}$ with $\vartheta^\prime, \vartheta^{\prime\prime} \in \Omega(X)$. The ``curvature'' of $\omega$ is the cochain $F = \beta \omega + \omega^2$, which is explicitly given by $F(\emptyset) = H = \slashed{\Dirac}^2 = \Dirac^2/2$, and
\begin{equation} \label{DefinitionF}
\begin{aligned}
  F(\vartheta) &~\defeq~ \cd(d \vartheta^\prime) - [\slashed{\Dirac}, \cd(\vartheta^\prime)] - \cd(\vartheta^{\prime\prime}), \\
  F(\vartheta_1, \vartheta_2) &~\defeq~ (-1)^{|\vartheta_1^\prime|}\bigl(\cd({\vartheta}_1^\prime \wedge {\vartheta}_2^\prime)- \cd({\vartheta}_1^\prime)\cd({\vartheta}_2^\prime)  \bigr),\\
  F(\vartheta_1, \dots, \vartheta_k) &~\defeq~ 0 \qquad \text{for} \quad k \geq 3.
  \end{aligned}
\end{equation}
In the first formula, the commutator is the {\em graded} commutator, keeping in mind that the Dirac operator is odd.
The corresponding \emph{Chern character} is then defined by the Chern-Weyl-type expression 
\begin{equation} \label{ChernExponential}
\Ch_{\Dirac} = \Str (e^{-F}),
\end{equation} 
where $\Str$ is the supertrace induced by \eqref{DefinitionSupertrace} on operators in $L^2(X, \Sigma)$ commuting with the right $\Cl_n$ action; see \eqref{OperatorSupertrace}. 

As the components of $F$ are unbounded, one needs to take some care making sense of the exponential $\Phi_t =  e^{-tF}$; an extensive discussion can be found in $\S4$ of \cite{GueneysuLudewig}.
For the purposes of this paper, it is convenient to define $\Phi_t$ by the recursive relation
\begin{equation} \label{PhiTRecursion}
\begin{aligned}
  \Phi_t(\vartheta_1, \dots, \vartheta_N) &= -\int_0^t e^{-(t-s)H} F(\vartheta_1) \Phi_s(\vartheta_2, \dots, \vartheta_N) \dd s \\
  &~~~
  - \int_0^t e^{-(t-s)H} F(\vartheta_1, \vartheta_2) \Phi_s(\vartheta_3, \dots, \vartheta_N) \dd s.
\end{aligned}
\end{equation}
for $N\geq 2$ and
\begin{equation*}
  \Phi_t(\emptyset) = e^{-tH}, \qquad \Phi_t(\vartheta_1) = - \int_0^t e^{-(t-s)H} F(\vartheta_1) e^{-sH} \dd s.
\end{equation*}
Since $e^{-\tau H}$ is a smoothing operator for every $\tau>0$ and the $F(\vartheta_1)$, $F(\vartheta_1, 2)$ are differential operators, an induction argument shows that also the operators  $\Phi_t(\vartheta_1, \dots, \vartheta_N)$ are smoothing operators whenever $t>0$. 
We can therefore take its supertrace componentwise, giving the components of the Chern character as
\begin{equation} \label{DefinitionChernCharacter}
\Ch_{\Dirac}(\vartheta_1, \dots, \vartheta_N) = \Str \Phi_1(\vartheta_1, \dots, \vartheta_N).
\end{equation}
It is not hard to work out that explicitly, $\Phi_t$ is given by the combinatorial formula
\begin{equation} \label{QuantizationMap}
  \Phi_t(\vartheta_1, \dots, \vartheta_N) \defeq \sum_{\substack{k=1 \\1 \leq a_1 < \dots < a_k\leq N}}^N  \!\!\!\!\!\! (-t)^k \int_{\Delta_k} \!\!\! e^{- t\tau_1 H} \prod_{i=1}^k F(\vartheta_{a_{i-1}+1}, \dots, \vartheta_{a_i}) e^{-t(\tau_{i+1}-\tau_{i}) H} \dd \tau,
\end{equation}
which is formula (4.8) of \cite{GueneysuLudewig}.
Using this formula, one checks using the cyclic permutation property of the supertrace that $\Ch_{\Dirac}$ is a {\em cyclic} cocycle, meaning that
\begin{equation} \label{ChernCharacterCyclic}
  \Ch_{\Dirac}(\vartheta_1, \dots, \vartheta_N) = (-1)^{n_k(n_N-n_k)} \Ch_{\Dirac}(\vartheta_{k+1}, \dots, \vartheta_N, \vartheta_1, \dots, \vartheta_k), 
\end{equation}
where $n_k = |\vartheta_1| + \dots + |\vartheta_k| - k$.

\begin{remark}
In \cite{GueneysuLudewig}, the Chern character is generally defined for Fredholm modules $\mathscr{M} = (\mathcal{H}, c, Q)$ over locally convex dg algebras $\Omega$, where $\mathcal{H}$ is a $\Z_2$-graded Hilbert space, $Q$ is an odd operator on $\mathcal{H}$ and $c: \Omega \rightarrow \mathcal{H}$ is a parity-preserving map (subject to several conditions). In the case that $n=\dim(X)$ is even, the above is the construction of of \cite{GueneysuLudewig}, applied to the Fredholm module given by $\mathcal{H} = L^2(X, \Sigma)$, $Q =\slashed{\Dirac}$ and $c(\vartheta) = \cd(\vartheta)$ (see in particular Example~2.3 of \cite{GueneysuLudewig}). The factors of $2^{-1/2}$ are present in order to adhere to the convention (standard in both stochastic analysis and physics) that the infinitesimal generator of Brownian motion is $\Delta/2$, and not $\Delta$.

In the odd-dimensional case, the trace considered above is not the standard supertrace on $\mathcal{B}(\mathcal{H})$, but still a tracial functional on the subalgebra $\mathcal{B}_{\Cl_n}(\mathcal{H}) \subset \mathcal{B}(\mathcal{H})$, and the constructions carry over essentially without changes.
\end{remark}

\paragraph{Properties of the Chern character.} The Chern character is coclosed, meaning that 
\begin{equation} \label{ChCoclosed}
  \Ch_{\Dirac}\bigl(b(c)\bigr) = 0 \qquad \text{for all} \qquad c \in \B^\natural(\Omega_\T(X)).
\end{equation}
This follows from \cite[Thm.~A]{GueneysuLudewig}; see also the short formal proof \cite[Prop.~1.6]{ludewig2020short}. Moreover, one can also show that  $\Ch_{\Dirac}$ extends by continuity to a linear functional on the entire complex $\B_\epsilon(\Omega_\T(X))$ \cite[Thm.~B]{GueneysuLudewig}, a fact that will be important when trying to plug in Bismut-Chern characters in \S\ref{SectionBismutChernCharacters}. 

The most remarkable property of $\Ch_{\Dirac}$, however, is that it satisfies the localization formula
\begin{equation} \label{LocalizationFormula}
  \Ch_{\Dirac}(c) = (2 \pi)^{-n/2} \int_X \widehat{A}(X) \wedge i(c)
\end{equation}
for any closed cocycle $c \in \B^\natural_\epsilon(\Omega_\T(X))$. Here $\widehat{A}(X)$ is the A-hat-genus form of $X$ (see e.g.\ \cite[formula (1.39)]{BGV}) and $i: \B(\Omega_\T(X)) \rightarrow \Omega(X)$ is the map defined by
\begin{equation*}
  i(\theta_1, \cdots, \theta_N) = \frac{(-1)^N}{N!} \,\theta_1^{\prime\prime} \wedge \cdots \wedge \theta_N^{\prime\prime}.
\end{equation*} 
Observe that we have $i = j^* \circ \rho$, where $\rho$ is the iterated integral map and $j: X \rightarrow \L X$ is the inclusion as constant loops.
A proof of \eqref{LocalizationFormula} can be found in \cite{ludewig2020short}; see also \cite[Thm.~E]{GueneysuLudewig}. 
There it is proved in the case that $\dim(X)$ is even, using the complex spinor bundle. 
This results in additional factors of the imaginary unit $i$, coming from the comparison formula \eqref{TraceComparisonEven}. 
When using the real spinor bundle, the proofs in the literature generalize in a straightforward way, also to the odd-dimensional case.
Observe that then the supertrace \eqref{DefinitionSupertrace} is odd, resulting in $\Ch_{\Dirac}$ being an odd functional. 
This matches the fact that also the right hand side of \eqref{LocalizationFormula} is an odd functional if $\dim(X)$ is odd.

\section{The path integral} \label{SectionPathIntegral}

In this section, we will define a \emph{path integral map}, which is a rigorous version of the formal definition
\begin{equation} \label{PathIntegralLater}
  \INT(\theta) \stackrel{\text{formally}}{=} \int_{\L X} e^{-S-\omega} \wedge \theta.
\end{equation}
This functional $\INT$ is then compared to the Chern character $\Ch_{\Dirac}$ introduced above.

\subsection{Definition of the path integral} \label{SectionDefinitionPathIntegral}

We now aim to define the path integral for all differential forms on $\L X$ that are (extended) iterated integrals, in other words, all forms in the image of the iterated integral map $\rho$ as defined in \eqref{ExtendedIteratedIntegralMap}.

\paragraph{A formal manipulation.}

A possible way to define the integral of a differential form $\theta$ on a finite-dimensional, oriented \emph{Riemannian} manifold $Y$ is the following: Take it to be the integral of its \emph{top degree component}, that is, the \emph{function} obtained by pairing $\theta$ with the volume form.

When trying to apply this to the loop space, we observe that while a Riemannian volume measure $\dd \gamma$ does not exist, it is well-known that the \emph{composite} measure $e^{-S} \dd \gamma$ does have a sensible interpretation as the Wiener measure (rigorously, this can be made sense of for example by using finite-dimensional approximation \cite{AnderssonDriver, BaerPfaeffle, MR3663619}). Formally rewriting \eqref{PathIntegralLater}, we therefore get
\begin{equation} \label{FormalManipulation}
\int_{\L X} e^{-S-\omega} \wedge \theta \stackrel{\text{formally}}{=} \int_{\L X} \qq(\theta) e^{-S}\dd \gamma \stackrel{\text{formally}}{=} \int_X \mathbb{E}_x \bigl[ \qq(\theta) ~\bigr|~ \xx_1 = x \bigr],
\end{equation}
where $\qq(\theta)$ is the (yet to be defined) ``top degree component'' of the composite form $e^{-\omega} \wedge \theta$, and on the right hand side, the conditional expectation is taken along a Brownian motion $\xx_\bullet$ in $X$ satisfying $\xx_0 = \xx_1 = x$. This leaves the task of finding an interpretation of this ``top degree component'' as an integrable function with respect to the Wiener measure.

\paragraph{The top degree component of iterated integrals.} 

 An interpretation of the ``top degree component'' on right hand side \eqref{FormalManipulation} has been found in our paper \cite{HanischLudewig2}, using formulas for such top degrees in finite dimensions and arguing by analogy. Our main result was then a reformulation of the resulting formula in terms of spin geometry. This requires the manifold to be spin, which is a natural condition as it is well-known to imply the ``orientability'' of the loop space \cite{StolzTeichner, Waldorf1}. Explicitly, we showed that for wedge products $\theta_1 \wedge \cdots \wedge \theta_N$ of one-forms $\theta_a \in L^2(\T, \gamma^* T^\prime X) \subset T_\gamma'\L X $, a sensible definition for their top degree component is
\begin{equation} \label{Ansatz}
 \qq(\theta_1 \wedge \cdots \wedge \theta_N) \defeq (-1)^N \sum_{\sigma \in S_N} \sgn(\sigma) \int_{\Delta_N}\str \left( [\gamma\|_{\tau_1}^0]^\Sigma \prod_{a=1}^N \cd\bigl(\theta_{\sigma_a}(\tau_a)\bigr) [\gamma\|_{\tau_{a+1}}^{\tau_a}]^\Sigma \right) \dd \tau,
\end{equation}
where $S_N$ is the $N$-th symmetric group and $\str$ is the supertrace \eqref{DefinitionSupertrace}. 

\begin{lemma} \label{Lemmabandq}
Interpreting the right hand side of \eqref{Ansatz} in the distributional sense, the functional $\qq$ can be applied pathwise to iterated integrals, giving the identity
\begin{equation} \label{TopDegreeVsQ}
  \qq\bigl(\rho(\vartheta_1 \wedge \cdots \wedge \vartheta_N)\bigr) = \str q(\vartheta_1, \dots, \vartheta_N),
\end{equation}
with $\vartheta_a \in \Omega_\T(X)$ and
\begin{equation} \label{ExplicitFormulaq}
 q(\vartheta_1, \dots, \vartheta_N) ~\defeq~ \int_{\Delta_N}  [\gamma\|_{\tau_1}^0]^\Sigma \prod_{a=1}^N \Bigl(\cd\bigl(\vartheta_a^{\prime\prime}(\gamma_{\tau_a})\bigr) - \cd\bigl(\iota_K \vartheta_a^\prime(\gamma_{\tau_a})\bigr) \Bigr) [\gamma\|_{\tau_{a+1}}^{\tau_a}]^\Sigma \dd \tau.
\end{equation}
\end{lemma}

\begin{proof}
Let $\theta \in \Omega^\ell(\L X)$. Then for each loop $\gamma$, $\theta_\gamma$ is an element of $L_{\mathrm{alt}}^\ell(TX, \R)$, which can be identified with the space $\mathscr{D}^\prime(\T^\ell, \gamma^*T^\prime X^{\boxtimes \ell})$  of distributions on $\T^\ell$ with values in the $\ell$-fold exterior product bundle $\gamma^*T^\prime X^{\boxtimes \ell}$ over $\T^\ell$ \cite[Lemma~1.2]{HanischLudewig2}. In particular, if $\vartheta \in \Omega^\ell(X)$, then the form $\vartheta(\tau) \in \Omega^\ell(\L X)$ defined in \eqref{TrivialLift} is pathwise a $\delta$-distribution supported on the point $(\tau, \dots, \tau) \in \T^\ell$.

More generally, if $\vartheta_a \in \Omega^{\ell_a+1}_\T(X)$, then for each $\tau_1 < \dots < \tau_N$, the integrand appearing in the formula \eqref{ExtendedIteratedIntegralMap} for the iterated integral $\rho(\vartheta_1, \dots, \vartheta_N)$ is a delta distribution supported at the point
\begin{equation*}
  (\underbrace{\tau_1, \dots, \tau_1}_{\ell_1}, \underbrace{\tau_2, \dots, \tau_2}_{\ell_2}, \dots, \underbrace{\tau_N, \dots, \tau_N}_{\ell_N}) \in \T^{\ell_1 + \dots + \ell_N}.
\end{equation*}
Explicitly, inserting these integrands into $\qq$ (interpreted in the distributional sense), we obtain
\begin{equation*}
 \qq\left( \bigwedge_{a=1}^N \bigl(\iota_K \vartheta_a^{\prime}(\tau_a) - \vartheta_a^{\prime\prime}(\tau_a)\bigr)\right) 
 =  \str \left( [\gamma\|_{\tau_1}^0]^\Sigma \prod_{a=1}^N \bigl( \cd\bigl(\vartheta_a^{\prime\prime}(\gamma_{\tau_a})\bigr)- \iota_K\cd\bigl(\vartheta^\prime_a(\gamma_{\tau_a}) [\gamma\|_{\tau_{a+1}}^{\tau_a}]^\Sigma \right),
\end{equation*}
where the factor of $(-1)^N$ has been pulled into the product, the integral disappears as we are ``integrating'' against a delta distribution and the sum over all permutations vanishes as the integrand is supported in $\Delta_N$ and hence all summands corresponding to a non-trivial $\sigma \in S_N$ are zero.
As $\rho(\vartheta_1, \dots, \vartheta_N)$ is the integral of these integrands over $\Delta_N$, the result follows.
\end{proof}

The discussion above is our motivation for taking the supertrace of the expression \eqref{ExplicitFormulaq} as an \emph{Ansatz} for the top degree functional of iterated integrals on a spin manifold $X$. 

We  would now like to integrate this function with respect to the Wiener measure. The obstacle to overcome here is that the Wiener measure is not defined on the smooth loop space $\L X$ but on the larger space $\L_c X$, to which \eqref{ExplicitFormulaq} does \emph{not} have a continuous extension. 
The problem here is two-fold: 
\begin{enumerate}
\item[(1)] The parallel transport appearing in \eqref{ExplicitFormulaq} is only defined for sufficiently regular paths, which form a zero set with respect to the Wiener measure. 
\item[(2)] The terms including $\iota_K\vartheta_a^\prime$ do not make sense along Brownian paths, since again, the vector field $K(\gamma) = \dot{\gamma}$ is only defined for sufficiently regular paths. 
\end{enumerate}
However, \eqref{ExplicitFormulaq} \emph{does} have a stochastic interpretation, where the parallel transport is interpreted in the stochastic sense, and the terms involving insertion of the velocity vector field are interpreted as Stratonovich stochastic integrals. In this sense, \eqref{ExplicitFormulaq} does have a \emph{measurable} extension, as we explain now.

\paragraph{Stochastic interpretation.} In the following, let $X$ be a compact spin manifold with spinor bundle $\Sigma$. Let moreover $\xx_\bullet$ be a Brownian motion in $X$, and denote by $[\xx\|_\bullet^\bullet]^\Sigma$ the \emph{stochastic} parallel transport in the spinor bundle along $\xx_\bullet$ (see e.g.\  \cite[\S2.3]{Hsu}, \cite[\S VIII]{Emery} or \cite[\S2.3]{GueneysuPhd}). In complete formal analogy to \eqref{ExplicitFormulaq}, given $\vartheta_1, \dots, \vartheta_N \in \Omega_\T(X)$, we consider the iterated Stratonovich integral
\begin{equation} \label{DefinitionQ}
  \tilde{q}(\vartheta_1, \dots, \vartheta_N) \defeq \int_{\Delta_N} [\xx_\bullet\|_{\tau_1}^0]^\Sigma \prod_{a=1}^N \Bigl( \cd\bigl(\vartheta_a^{\prime\prime}(\xx_{\tau_a})\bigr) \dd \tau_a - \cd\bigl(\iota_\bullet \vartheta_a^\prime(\xx_{\tau_a})\bigr) * \dd \xx_{\tau_a}\Bigr) [\xx_\bullet\|_{\tau_{a+1}}^{\tau_a}]^\Sigma,
\end{equation}
adapted to the process $\xx_\bullet$. Here and throughout, we write $*\dd \xx$ for Stratonovich integration and $\dd \xx$ for It\^o integration.
 
In fact, this integral is the \emph{canonical} stochastic extension of \eqref{ExplicitFormulaq}, as we explain now. 
Observe first that \eqref{ExplicitFormulaq} is welldefined also if the loop $\gamma$ is only piecewise smooth. Now for a subdivision $\tau = \{0 = \tau_0 < \tau_1 < \dots < \tau_m = 1\}$ of the interval $[0, 1]$, let $\xx_\bullet^\tau$ be a stochastic process such that $\xx_{\tau_j}^\tau = \xx_{\tau_j}$ almost surely for each $j=1, \dots, m$ and such that almost surely, the sample paths of $\xx_\bullet^\tau$ are minimizing geodesics on each of the subintervals $[\tau_{j-1}, \tau_j]$. 
It then follows from the finite approximation property of  the stochastic parallel transport and the Stratonovich integral (see \cite[7.14 \& 8.15]{Emery}) that we have
\begin{equation} \label{QApproximation}
 \lim_{|\tau|\rightarrow 0} q(\vartheta_1, \dots, \vartheta_N)_{\xx^\tau_\bullet} = \tilde{q}(\vartheta_1, \dots, \vartheta_N)
\end{equation}
in probability, where the left hand side denotes the pathwise evaluation of the $q$-functional on the piecewise smooth sample paths of $\xx_\bullet^\tau$, and the limit is taken over any sequence of partitions $\tau$ the mesh of which tends to zero. 

\begin{lemma} \label{LemmaIndependentOfRepresentation}
  The random variable $\tilde{q}(\vartheta, \dots, \vartheta_N)$ depends only on the differential form $\rho(\vartheta, \dots, \vartheta_N) \in \Omega(\L X)$ and not on the particular representation as iterated integral.
\end{lemma}

\begin{proof}
It follows directly from \eqref{TopDegreeVsQ} that $\str q(\vartheta_1, \dots, \vartheta_N)$ only depends on  $\rho(\vartheta_1, \dots, \vartheta_N)$ and not on the particular choice of $\vartheta_1, \dots, \vartheta_N$. One can now define a generalization of the functional $\qq$ defined at general (non-looping) paths $\gamma$, by dropping the supertrace in \eqref{Ansatz}; it is then straightforward to check that the result and the proof of Lemma~\ref{Lemmabandq} extend to this more general context. We obtain that $q$ does only depend on the iterated integral in $\Omega(\L X)$, not the particular representation. By \eqref{QApproximation}, the same must be true for  $\tilde{q}$.
\end{proof}

By the above lemma, for every $\theta = \rho(c)$, $c \in \mathsf{B}^\natural(\Omega_{\T}(X))$, we obtain a well-defined random variable $\tilde{\qq}(\theta)$ adapted to the Brownian bridge $\xx_\bullet$ with $\xx_0 = \xx_1$. Explicitly, it is defined by
\begin{equation*}
  \tilde{\qq}(\theta) = \str \tilde{q}(c); 
\end{equation*}
by Lemma~\ref{LemmaIndependentOfRepresentation}, the definition is independent of the choice of $c$. 
By the estimate from Lemma~\ref{LemmaEstimateQt} below (applied to the special case \eqref{SpecialCaseOmegaV}), $\qq(\theta)$ is indeed integrable whenever $\theta = \rho(\vartheta_1, \dots, \vartheta_N)$ for some $\vartheta_1, \dots, \vartheta_N$. By definition \eqref{EntireSeminorms} of the entire seminorm, this estimate shows moreover that $\qq(\theta)$ extends to a well-defined integrable process for all $\theta \in \Omega_{\mathrm{int}}(LX)$.

\begin{definition}[The path integral] \label{DefinitionPathIntegral}
The \emph{path integral map}
\begin{equation*}
  \INT : \Omega_{\mathrm{int}}(\L X) \longrightarrow \R,
\end{equation*}
is defined by the formula
\begin{equation} \label{StochasticPathIntegral}
  \INT\bigl(\theta\bigr) \defeq \int_X \mathbb{E}_x \left[ \tilde{\qq}(\theta) \cdot \exp \left( -\frac{1}{8} \int_0^1 \mathrm{scal}(\xx_t) \dd t \right) ~\Bigl|~ \xx_1 = x\right] \dd x,
\end{equation}
where we take the conditional expectation of a Brownian motion $\xx_\bullet$ starting at $x$, constrained to end at the same point.
\end{definition}

\begin{remark}$-$ \label{RemarkScalarCurvature} 
\begin{enumerate}
\item[(1)]
Our definition \eqref{StochasticPathIntegral} differs from the heuristic version in the right hand side of \eqref{FormalManipulation} by a scalar curvature term.
This scalar curvature term is a necessary addition to the formula in order to ensure supersymmetry. In effect, it comes from the Lichnerowicz formula 
\begin{equation} \label{LichnerowiczFormula}
  \slashed{\Dirac}^2 = \frac{1}{2}\nabla^*\nabla + \frac{1}{8} \mathrm{scal}.
\end{equation}
\item[(2)] Formula \eqref{StochasticPathIntegral} provides (up to factors) an extension of the path integral defined by Lott \cite[\S V]{Lott} to a significantly larger domain. There, the scalar curvature term is stated to arise due to ``quantum effects'', coming from factor ordering in the Hamiltonian form of the path integral and the Lichnerowicz formula (see ibid.\ 634).
\end{enumerate}
\end{remark}

\subsection{A Feynman-Kac formula} \label{SectionFeynmanKac}

In this section, we will prove a  Feynman-Kac type formula that we need in order to compare the path integral map $\INT$ defined in \S\ref{SectionDefinitionPathIntegral} with the Chern character $\Ch_{\Dirac}$. This formula allows us to express the kernels of certain operators on $L^2(X, \Sigma)$ as expectation values of corresponding stochastic processes, where $X$ is a compact spin manifold with spinor bundle $\Sigma$. We remark that the formulas we obtain here have generalizations for operators acting on sections of general vector bundles over complete but not necessarily compact Riemannian manifolds; see \cite{BoldtGueneysu}.

Throughout, fix a collection of 1-forms $\omega_1, \dots, \omega_N \in \Omega^1(X, \End_{\Cl_n}(\Sigma))$. Moreover, we let $\omega_a^\sharp \in C^\infty(X, TX \otimes \End_{\Cl_n}(\Sigma))$, and $\div \omega_a \in C^\infty(X, \End_{\Cl_n}(\Sigma))$ be defined by
\begin{equation*}
 \omega_a^\sharp \defeq \sum_{j=1}^n e_j \otimes \omega_a[e_j], \qquad \div \omega_a(x) \defeq \sum_{j=1}^n (\nabla_{e_j}^\Sigma \omega_a) [e_j].
\end{equation*}
These definitions are independent of the choice of the orthonormal basis $e_1, \dots, e_n$ of $T_xX$. Observe that for $x \in X$, $\omega^\sharp(x)$ may be interpreted as a linear map $T^*_xX \otimes \Sigma_x \rightarrow \Sigma_x$.
Moreover, fixing a collection of \emph{potentials}, by which we mean sections $V_1, \dots, V_N \in C^\infty(X, \End_{\Cl_n}(\Sigma))$, we define a stochastic process $Q_\bullet$ by the iterated Stratonovich integral
\begin{equation} \label{DefinitionGeneralQ}
  Q_t = \int_0^t \int_0^{\tau_N} \cdots \int_0^{\tau_2} [\xx\|_{\tau_1}^0]^\Sigma \prod_{a=1}^N \Bigl( \omega_a(\xx_{\tau_a}) * \dd \xx_{\tau_a} + V_a(\xx_{\tau_a}) \dd \tau_a\Bigr) [\xx\|_{\tau_{a+1}}^{\tau_a}]^{\Sigma}.
\end{equation}
The following lemma follows from writing $Q_\bullet$ as an It\^o integral using the Stratonovich-to-It\^o formula (see Lemma~\ref{LemmaQ} below) and the It\^o isometry.

\begin{lemma} \label{LemmaEstimateQt}
For each $t>0$, there exists $C>0$ such that
\begin{equation*}
  \mathbb{E}_x\bigl[|Q_t|^2\bigr] \leq \frac{C^N}{\lfloor N/2\rfloor!^2} \prod_{a=1}^N \bigl( \|V_a\|_\infty^2 + \|\omega_a\|_\infty^2 + \|\div \omega_a\|_\infty^2 \bigr),
\end{equation*}
where the constant $C>0$ is independent of the choice of $\omega_a$ and $V_a$ and the expectation value is taken with respect to a Brownian motion starting at $x \in X$.
\end{lemma}

By this lemma, if $\xx_t$ is a Brownian motion starting at $x \in X$, then for any $\psi \in L^2(X, \Sigma)$ and any $t\geq0$, the $\Sigma_x$-valued random variable $S_t \cdot Q_t \cdot \psi(\xx_t)$ has a finite expectation value, where the \emph{scalar curvature process} $S_\bullet$ is defined by
\begin{equation} \label{ScalarCurvatureProcess}
S_t = \exp \left( -\frac{1}{8} \int_0^t \mathrm{scal}(\xx_s) \dd s \right).
\end{equation}
To give an operator theoretic formula for this expectation values, define the auxiliary operators
\begin{equation} \label{OtherFs}
A_a = \omega_a^\sharp \nabla + \frac{1}{2} \div \omega_a + V_a, \qquad A_{a, a+1} = \frac{1}{2} \omega_a^\sharp \omega_{a+1},\qquad A_{a, \dots, a+m} = 0, \quad m\geq 2.
\end{equation}
Using these, inductively define operators $\Psi_t^{(a)}$ by setting $\Psi_t^{(N)} = e^{-tH}$, $\Psi^{(N+1)}_t = 0$ and
\begin{equation} \label{RecursiveGeneralPhiT}
\Psi_t^{(a)} = \int_0^t e^{-(t-s)H} A_{a+1} \Psi_s^{(a+1)} \dd s + \int_0^t e^{-(t-s)H} A_{a+1, a+2} \Psi_s^{(a+2)} \dd s
\end{equation} 
for $a < N$, where $H = \frac{1}{2} \Dirac^2$. Observe that since each of the operators defined in \eqref{OtherFs} is a differential operator (of order at most one) and $e^{-tH}$ is smoothing, each of the $\Psi_t^{(a)}$ ends up being a smoothing operator as the composition of a smoothing operator and a differential operator is again smoothing. It will be important that the operators $\Psi_t^{(a)}$ solve the differential equation
\begin{equation} \label{OperatorODE}
 \tfrac{\dd}{\dd t}\Psi_t^{(a)} = - H \Psi_t^{(a)} + A_{a+1} \Psi_t^{(a+1)} + A_{a+1, a+2} \Psi_t^{(a+2)}, \qquad \Psi_0^{(a)} = \begin{cases} 0 & a \neq N \\ \id & a = N \end{cases},
 \end{equation}
 which can be taken as an alternative definition for the $\Psi_t^{(a)}$.

\begin{proposition}[Feynman-Kac formula] \label{PropFeynmanKac}
For any $\psi \in L^2(X, \Sigma)$, any $x \in X$ and any $t > 0$, we have
\begin{equation*}
\mathbb{E}_x\bigl[S_t \cdot Q_t \cdot \psi(\xx_t) \bigr] =  (\Psi_t^{(0)} \psi)(x),
\end{equation*}
where on the left hand side, we take the expectation value of a Brownian motion $\xx_\bullet$ such that $\xx_0 = x$.
\end{proposition}

We remark that in the case that $N=1$, a similar result for not necessarily compact manifolds and arbitrary bundles was recently obtained in \cite{BoldtGueneysu}.

\medskip

In order to prove Prop.~\ref{PropFeynmanKac}, it is convenient to lift the process $Q_t$ to the $\mathrm{Spin}_n$-principal bundle $P {\rightarrow} X$ provided by the spin structure. The advantage of this approach is that we can pass from sections of bundles over $X$ to functions with values in a fixed vector space. 
More precisely, the elements $V_a$ and $\omega_a$ correspond to elements $\bfV_a \in C^\infty(P, \Cl_n)^{\Spin_n}$ and $\bfomega \in \Omega^1(P, \Cl_n)^{\Spin_n}$, defined by
\begin{equation*}
  V_a(x) = [u, \bfV_a(u)], \qquad \omega_a(x) = [u, \bfomega_a(u)]
\end{equation*}
for any $u \in P_x$. Here we represented $\End_{\Cl_n}(\Sigma) \cong \Cl(TX) $ as the associated bundle $\End_{\Cl_n}(\Sigma) = P \times_{\Spin_n} \Cl_n$ (where in this case, $\Spin_n$ acts on $\Cl_n$ through $\mathrm{SO}_n$). Alternatively, we can view elements of $P_x$ as isometries $u: \Cl_n \rightarrow \Sigma_x$ respecting the right $\Cl_n$-action, which leads to the relations
\begin{equation} \label{RelationsVandOmega}
  V_a(x) = u \, \bfV_a(u) \, u^*, \qquad \omega_a(x) w = u \bigl( \bfomega(u) w^{\mathrm{hor}} \bigr) u^*,
\end{equation}
with $w^{\mathrm{hor}}$ the horizontal lift of $w \in T_xX$.
Let $\UU_\bullet$ be the horizontal lift to $P$ of a Brownian motion $\xx_\bullet$ (this can be defined as the solution to a stochastic differential equation, as discussed in \cite[\S2.3]{Hsu}; see also \cite[\S2.2]{GueneysuPhd}). 
Define the $\Cl_n$-valued stochastic processes $\bfQ_\bullet^{(a)}$, $a=1, \dots, N$, inductively as the solution to the Stratonovich differential equation
\begin{equation}
  \dd \bfQ^{(a)}_t = \bfQ^{(a-1)}_t \cdot \Bigl( \bfomega_a(\UU_t) * \dd \UU_t + \bfV_a(\UU_t) \dd t \Bigr), \qquad \bfQ_0^{(a)} = 0,
\end{equation}
where we let $\bfQ^{(0)}_t \equiv 1$ and (for use in the lemma below) $\bfQ^{(-1)}_t \equiv 0$. 
Explicitly, the solution for $a > 0$ to this differential equation is
\begin{equation*}
\bfQ^{(a)}_t = \int_0^t \int_0^{\tau_{a}} \cdots \int_0^{\tau_{2}} \Bigl( \bfomega_1(\UU_{\tau_1}) * \dd \UU_{\tau_1} + \bfV_1(\UU_{\tau_1}) \dd \tau_1 \Bigr) \cdots \Bigl( \bfomega_a(\UU_{\tau_a}) * \dd \UU_{\tau_a} + \bfV_a(\UU_{\tau_a}) \dd {\tau_a} \Bigr).
\end{equation*}
In terms of $\UU_\bullet$, the stochastic parallel transport is just given by $[\xx\|_s^t]^\Sigma = \UU_t \UU_s^*$, so it follows with a view on \eqref{RelationsVandOmega} that we have
\begin{equation} \label{FatQvsQ}
\UU_0 \bfQ_t^{(N)} \UU_t^*  =  Q_t,
\end{equation} 
the process defined in \eqref{DefinitionGeneralQ}. 

\begin{lemma} \label{LemmaQ} 
  The processes $\bfQ_\bullet^{(a)}$ solve the It\^o differential equation
  \begin{equation} \label{StratonovichODEE}
   \dd \bfQ^{(a)}_t = \bfQ^{(a-1)}_t \cdot \Bigl( \bfomega_a(\UU_t) \dd \UU_t + \bfV^\prime_a(\UU_t) \dd t \Bigr) + \frac{1}{2}\bfQ^{(a-2)}_t \cdot \bfomega_{a-1}^\sharp(\UU_t) \bfomega_a(\UU_t) \dd t,
  \end{equation}
where $\bfV^\prime_a  \in C^\infty(P, \Cl_n)^{\Spin_n}$ is the lift to $P$ of the section $V_a^\prime := V_a + \frac{1}{2} \div \omega_a$.
\end{lemma}

\begin{remark}
  This lemma can be used to transfer the multiple Stratonovich integral $Q_t$ given by \eqref{DefinitionGeneralQ} into a multiple It\^o integral. After a close look at the combinatorics involved, one observes that the result is a sum of several integrals over various lower-dimensional simplices, formally similar to \eqref{QuantizationMap}. This is a general version of so-called \emph{Hu-Meyer formulas} \cite{HuMeyer1}.
\end{remark}

\begin{proof}
First let $a \geq 2$.
Then by the Stratonovich-to-It\^o formula $X * \dd Y = X \dd Y + \frac{1}{2}\dd X \dd Y$ and the Ito formula, we have
\begin{equation*}
\begin{aligned}
  2 \bigl(\bfQ_t^{(a-1)} \cdot \bfomega_a&(\UU_t) * \dd \UU_t - \bfQ_t^{(a-1)} \cdot \bfomega_a(\UU_t) \dd \UU_t\bigr) = \dd\bigl[ \bfQ_t^{(a-1)} \cdot \bfomega_a(\UU_t) \bigr] \dd \UU_t \\
  &= \dd \bfQ_t^{(a-1)} \cdot \bfomega_a(\UU_t) \dd \UU_t + \bfQ_t^{(a-1)} \cdot \dd[\bfomega_a(\UU_t)] \dd \UU_t + \underbrace{\dd \bfQ_t^{(a-1)} \dd[\bfomega_a(\UU_t)] \dd \UU_t}_{=0}\\
  &= \bfQ_t^{(a-2)} \cdot \bfomega_{a-1}(\UU_t) \dd \UU_t \cdot \bfomega_a(\UU_t) \dd \UU_t + \bfQ_t^{(a-1)} \cdot \nabla_{\dd \UU_t} \bfomega_a(\UU_t) \dd \UU_t\\
  &= \bfQ_t^{(a-2)} \cdot \bfomega_{a-1}^\sharp(\UU_t)\bfomega_a(\UU_t) \dd t + \bfQ_t^{(a-1)} \cdot \div \bfomega_a(\UU_t) \dd t,
\end{aligned}
\end{equation*}
which gives \eqref{StratonovichODEE}.
Here we use that the higher terms in the It\^o formula do not contribute as $\dd t \cdot \dd t = \dd t \cdot \dd \UU_t = 0$ and the formula
\begin{equation} \label{BilinearFormProcess}
  B[\dd \UU_t, \dd \UU_t] = \tr^{\mathrm{hor}} B(\UU_t) \dd t
\end{equation}
which holds for any symmetric 2-tensor $B$ on $P$; here the \emph{horizontal trace} is defined by
\begin{equation*}
  \tr^{\mathrm{hor}} B(u) = \sum_{j=1}^n B(E_j, E_j),
\end{equation*}
where $E_1, \dots, E_n$ is an orthonormal basis of the horizontal tangent space $T^{\mathrm{hor}}_u P$ at $u$.
Finally, for $a=1$, the same calculation works, except that $\bfQ^{(a-1)} \equiv 1$ and consequently $\dd \bfQ_t^{(a-1)} \equiv 0$. Since by convention, $\bfQ^{(a-2)} \equiv 0$, formula \eqref{StratonovichODEE} holds also in that case.
\end{proof}

Fix $\psi \in L^2(X, \Sigma)$ and $t > 0$. For $0 \leq s \leq t$ and $a = 1, \dots, N$, let $\psi_s^{(a)}$ be the section of $\Sigma$ defined by
\begin{equation} \label{DefinitionSmallPsi}
  \psi_s^{(a)}(x) = (\Psi_{t-s}^{(a)} \psi)(x).
\end{equation}
Observe that $\psi_s^{(a)}$ is smooth for $s \in [0, t)$.

\begin{lemma} \label{LemmaPsi}
Let $\bfpsi_s^{(a)} \in C^\infty(P, \Cl_n)^{\Spin_n}$ be the horizontal lift of $\psi_s^{(a)}$. Then the process $\bfpsi^{(a)}_s(\UU_s)$ satisfies
\begin{equation*}
 \dd [\bfpsi^{(a)}_s(\UU_s)] = d\bfpsi_s^{(a)} (\UU_s) \dd \UU_s 
  + \left(\frac{1}{8} \mathbf{scal}(\UU_s) \bfpsi^{(a)}_s - A_{a+1}^{\mathrm{hor}}\bfpsi_s^{(a+1)}(\UU_s) - A_{a+1, a+2}^{\mathrm{hor}} \bfpsi_s^{(a+2)}(\UU_s) \right)\dd s,
\end{equation*}
where $\mathbf{scal}$ is the lift of the scalar curvature function to $P$ and $A_a^{\mathrm{hor}}$, $A_{a-1, a}^{\mathrm{hor}}$ are the horizontal lifts of the operators $A_a$ and $A_{a-1, a}$.
\end{lemma}

Here, if $A$ is a differential operator on $X$ acting on sections of $\Sigma$, its horizontal lift (to $P$) is the operator $A^{\mathrm{hor}}$ on $C^\infty(P, \Cl_n)^{\mathrm{Spin}_n}$ such that $A^{\mathrm{hor}}\boldsymbol{\psi}$ is the lift to $P$ of $A \psi$, if $\boldsymbol{\psi}$ is the lift of $\psi$.

\begin{proof}
By \eqref{OperatorODE}, the functions $\bfpsi^{(a)}_s$ satisfy the system of differential equations
\begin{equation*}
\tfrac{d}{d s}\bfpsi^{(a)}_s = H^{\mathrm{hor}} \bfpsi_s^{(a)} - A^{\mathrm{hor}}_{a+1} \bfpsi_s^{(a+1)} - A^{\mathrm{hor}}_{a+1, a+2} \bfpsi_s^{(a+2)}.
\end{equation*}
Let $\boldsymbol{\Psi}$ be the function on $P \times [0, t]$ defined by $\boldsymbol{\Psi}(u, s) = \bfpsi^{(a)}_s(u)$ and let ${\mathsf{U}}_s = (\UU_s, s)$. Then $\boldsymbol{\Psi}(\mathsf{U}_s) = \bfpsi^{(a)}_s(\UU_s)$ and hence by the Stratonovich product rule, 
\begin{equation*}
\begin{aligned}
\dd [\bfpsi_s(\UU_s)] &= \dd [\boldsymbol{\Psi}(\mathsf{U}_s)] = d \boldsymbol{\Psi}(\mathsf{U}_s) * \dd \mathsf{U}_s = d \bfpsi^{(a)}_s(\UU_s) * \dd \UU_s + \tfrac{d}{d s} \bfpsi^{(a)}_s (\UU_s)\dd s \\
&= d \bfpsi^{(a)}_s(\UU_s) * \dd \UU_s + H^{\mathrm{hor}} \bfpsi^{(a)}_s(\UU_s) - A^{\mathrm{hor}}_{a+1} \bfpsi_s^{(a+1)}(\UU_s) - A^{\mathrm{hor}}_{a+1, a+2} \bfpsi_s^{(a+2)}(\UU_s).
\end{aligned}
\end{equation*}
On the other hand, the Stratonovich-to-It\^o formula yields
\begin{equation*}
  d \bfpsi^{(a)}_s(\UU_s) * \dd \UU_s = d \bfpsi^{(a)}_s(\UU_s) \dd \UU_s - \tfrac{1}{2} (\nabla^*\nabla)^{\mathrm{hor}} \bfpsi^{(a)}_s(\UU_s),
\end{equation*}
where $(\nabla^*\nabla)^{\mathrm{hor}}$ is the horizontal lift of the connection Laplacian on $\Sigma$. Combining this with the Lichnerowicz formula \eqref{LichnerowiczFormula} (keeping in mind that $H = \Dirac^2/2$) results in the claimed scalar curvature term; this finishes the proof.
\end{proof}

\begin{proof}[of Prop.~\ref{PropFeynmanKac}]
Define $\bfN_s^{(a)} := \bfQ^{(a)}_s \cdot \bfpsi_s^{(a)}(\UU_s)$. 
By the It\^o formula, 
\begin{equation*}
  \dd \bfN_s^{(a)} = \dd \bfQ_s^{(a)} \cdot \bfpsi_s^{(a)}(\UU_s) + \bfQ_s^{(a)} \cdot \dd\bigl[\bfpsi_s^{(a)}(\UU_s)\bigr] + \dd \bfQ_s^{(a)} \cdot \dd\bigl[\bfpsi_s^{(a)}(\UU_s)\bigr].
\end{equation*}
Each of the terms is readily calculated using Lemmas~\ref{LemmaQ} \& \ref{LemmaPsi}. In particular, the third term is
\begin{equation*}
  \bfQ_s^{(a-1)} \cdot \bfomega_a(\UU_s) \dd \UU_s \cdot d \bfpsi^{(a)}_s(\UU_s) \dd \UU_s = \bfQ_s^{(a-1)} \cdot \bfomega^\sharp(\UU_s) d \bfpsi^{(a)}_s(\UU_s) \dd s,
\end{equation*}
where we used the rule $\dd \UU_s \cdot \dd s = \dd s \cdot \dd s = 0$ and the identity \eqref{BilinearFormProcess} again. Collecting all the terms, we obtain
\begin{equation} \label{DifferentialNa}
\begin{aligned}
  \dd \bfN_s^{(a)} &= \bfQ^{(a-1)}_s \cdot \Bigl(\bfomega^\sharp d \bfpsi^{(a)}_s(\UU_s) + \bfV^\prime_a(\UU_s)   \bfpsi^{(a)}_s(\UU_s)\Bigr) \dd s - \bfQ_s^{(a)} \cdot A_{a+1}^{\mathrm{hor}} \bfpsi^{(a+1)}_s(\UU_s) \dd s \\
 & + \frac{1}{2}\bfQ^{(a-2)} \cdot (\bfomega_{a-1}^\sharp \bfomega_a)(\UU_s) \cdot \bfpsi_s^{(a)} \dd s
   - \bfQ^{(a)} \cdot  A_{a+1, a+2}^{\mathrm{hor}} \cdot \bfpsi_s^{(a+2)} \dd s\\
   & + \frac{1}{8} \bfQ^{(a)}_s \cdot \mathbf{scal}(\UU_s) \cdot \bfpsi_s^{(a)} \dd s + \bfQ_s^{(a-1)} \cdot \bfomega_a(\UU_s) \dd \UU_s \cdot \bfpsi_s^{(a)} + \bfQ_s^{(a)} \cdot d\bfpsi_s^{(a)} (\UU_s) \dd \UU_s.
\end{aligned}
\end{equation}
Now set 
\begin{equation*}
  \bfN_s \defeq  \mathbf{S}_s \cdot \sum_{a=0}^N  \bfN_s^{(a)}
\end{equation*}
where $\mathbf{S}_\bullet$ is the lift to $P$ of the scalar curvature process \eqref{ScalarCurvatureProcess}. We obtain
\begin{equation} \label{DifferentialNs1}
 \dd \bfN_s = \dd \mathbf{S}_s \cdot \sum_{a=0}^N \bfN_s^{(a)} + \mathbf{S}_s \cdot \sum_{a=0}^N \dd \bfN_s^{(a)} + \sum_{a=0}^N \dd \mathbf{S}_s \cdot \dd \bfN_s^{(a)}.
\end{equation}
Since $\mathbf{S}_\bullet$ satisfies the differential equation
\begin{equation*}
  \dd \mathbf{S}_s = - \frac{1}{8} \mathbf{S}_s \cdot \mathbf{scal}(\UU_s) \dd s,
\end{equation*}
the third term of \eqref{DifferentialNs1} vanishes (again using $\dd s \cdot \dd s = \dd s \cdot \dd \UU_s = 0$). The middle term can be calculated using \eqref{DifferentialNa}: In view of the definition \eqref{OtherFs} of the operators $A_a$, $A_{a, a+1}$, the terms in the first two rows of \eqref{DifferentialNa} telescope, while the scalar curvature terms cancel with the terms coming from the first summand of \eqref{DifferentialNs1}. 
We are left with
\begin{equation*}
  \dd \bfN_s = \sum_{a=0}^N \Bigl(\bfQ_s^{(a-1)} \cdot \bfomega_a(\UU_t) \dd \UU_t \cdot \bfpsi_s^{(a)}(\UU_s) + \bfQ_s^{(a)} \cdot d\bfpsi_s^{(a)} (\UU_s) \dd \UU_s \Bigr),
\end{equation*}
so that $\bfN_s$ is a Martingale by the Martingale property of the It\^o integral. Let 
\begin{equation*}
  N_s \defeq \UU_0 \bfN_s \UU_s^* = S_s \cdot \sum_{a=0}^N  \UU_0 \bfQ_s^{(a)} \UU_s^* \cdot \UU_s \bfpsi_s^{(a)}(\xx_s)\UU_s^* = S_s \cdot \sum_{a=0}^N Q_s^{(a)} \cdot \psi_s^{(a)}(\xx_s),
\end{equation*}
where we used \eqref{FatQvsQ} and \eqref{RelationsVandOmega}. By definition, $\psi_s^{(a)}(\xx_s) = (\Psi_{t-s}^{(a)}\psi)(\xx_s)$. We have
\begin{equation*}
  Q_0^{(a)} = \begin{cases} 1 & a = 0 \\ 0 & a \geq 1 \end{cases} 
  \qquad  \text{and} \qquad
  \psi_{t}^{(a)}(\xx_0) =  \begin{cases} \psi(x) & a = N \\ 0 & a < N, \end{cases}
\end{equation*}
using that $\xx_0 = x$. Therefore
\begin{equation*}
N_0  =  (\Psi^{(0)}_t \psi)(x)  \qquad \text{and} \qquad  N_t = {S}_t \cdot Q^{(N)}_t \cdot \psi(\xx_t),
\end{equation*}
and we obtain
\begin{equation*}
\mathbb{E}_x\bigl[ S_t \cdot Q^{(N)}_t \cdot \psi(\xx_t)\bigr] - (\Psi_t^{(0)} \psi)(x) 
  = \mathbb{E}_x\Bigl[ N_t - N_0 \Bigr] = \mathbb{E}_x\bigl[ \UU_0 \bfN_t \UU_t^*- \UU_0\bfN_0\UU_0^*\bigr] = 0,
\end{equation*}
as the expectation value of a Martingale is zero.
\end{proof}

\subsection{Comparison of path integral and Chern character} \label{SectionComparison}

\begin{theorem} \label{TheoremEqualityChPI}
Let $\vartheta_1, \dots, \vartheta_N \in \Omega_\T(X)$ and as always, write $\vartheta_a = \vartheta_a^{\prime} + dt \wedge \vartheta_a^{\prime\prime}$ with $\vartheta_a^\prime, \vartheta_a^{\prime\prime} \in \Omega(X)$. Then if for each pair of subsequent indices $a$, $a+1$, at least one of $\vartheta_a^\prime$ or $\vartheta_{a+1}^\prime$ has degree less or equal to one, then
\begin{equation} \label{EqualityChIThm}
  \Ch_{\Dirac}(\vartheta_1, \dots, \vartheta_N) = \INT\bigl(\rho(\vartheta_1, \dots, \vartheta_N)\bigr).
\end{equation}
\end{theorem}

Incidentally,  the Bismut-Chern characters defined in \S\ref{SectionBismutChernCharacters} below satisfy the assumptions of the theorem. 

\begin{remark}
From inspection of the arguments below, it is easy to see that in general, the difference of $\Ch_{\Dirac}(\vartheta_1, \dots, \vartheta_N)$ and $\INT\bigl(\rho(\vartheta_1, \dots, \vartheta_N)\bigr)$ can be written as an expectation value of iterated Stratonovich integrals depending on the differences 
\begin{equation*}
F(\vartheta_a, \vartheta_{a+1}) - \frac{1}{2}\sum_{j=1}^n \cd(\iota_{e_j} \vartheta_a) \cd(\iota_{e_j} \vartheta_{a+1}),
\end{equation*}
which happen to vanish if at least one of $\vartheta_a^\prime$ and $\vartheta_{a+1}^\prime$ has degree less or equal to one. 
This suggests that the fact that the equality \eqref{EqualityChIThm} breaks down in the case that these terms are non-zero comes from bad interactions of multiple stochastic integrals with the diagonal singularities of the top degree current \eqref{Ansatz}.
\end{remark}

A crucial ingredient of the proof is the following rather magical lemma, which connects the operators \eqref{OtherFs} to the terms appearing in the definition of $\Ch_{\Dirac}$. 

\begin{lemma} \label{ComparisonLemma}
For all $\vartheta \in \Omega(X)$, we have
  \begin{equation} \label{Comparison1}
    \cd( d\vartheta) - [\slashed{\Dirac}, \cd(\vartheta)] = \sum_{i=1}^n \cd(\iota_{e_i} \vartheta) \nabla_{e_i} - \frac{1}{2}\cd(d^* \vartheta),
  \end{equation}
where $d^*$ is the adjoint of $d$. Moreover, for $\vartheta_1, \vartheta_2 \in \Omega(X)$, if either $\vartheta_1$ or $\vartheta_2$ has degree at most one, then
  \begin{equation} \label{Comparison2}
    (-1)^{|\vartheta_1|}\bigl( \cd(\vartheta_1)\cd(\vartheta_2) - \cd(\vartheta_1 \wedge \vartheta_2)\bigr) = \frac{1}{2}\sum_{j=1}^n \cd(\iota_{e_j} \vartheta_1) \cd(\iota_{e_j} \vartheta_2).
  \end{equation}
\end{lemma}

\begin{proof}
If $|\vartheta|=0$, the right hand side of \eqref{Comparison1} vanishes, so this is just the usual commutator relation of the Dirac operator for functions; the case $|\vartheta_1|=1$ is Prop.~3.45 of \cite{BGV}. In the general case, the identity \eqref{Comparison1} can be found in \cite{BoldtGueneysu}.

Equation \eqref{Comparison2} is a pointwise statement. Hence we may choose an orthonormal basis $e_1, \dots, e_n$ of $T_x X$ and (by linearity) assume that $\vartheta_1 = e_i$ and $\vartheta_2 = e_{j_1} \wedge \cdots \wedge e_{j_\ell}$, where $J = \{1 \leq j_1 < \dots < j_\ell \leq n\}$ is an arbitrary multi-index. Observe that if $i \notin J$, then both sides of \eqref{Comparison2} are zero. If $i = j_a \in J$, then $\vartheta_1 \wedge \vartheta_2 = 0$ and
\begin{equation*}
 \cc(\vartheta_1) \cc(\vartheta_2) =   (-1)^a \cc(e_{j_1}) \cdots \widehat{\cc(e_{j_a})} \cdots \cc(e_{j_\ell}) = - \sum_{k=1}^n \cc(\iota_{e_k} e_i) \cc(\iota_{e_k} e_{j_1} \wedge \cdots \wedge e_{j_\ell}),
\end{equation*}
where we used that all summands on the right hand side vanish except that with $k=a$.
Swapping the roles of $\vartheta_1$ and $\vartheta_2$, the right hand side of \eqref{Comparison2} remains unhanged, while the left hand side picks up a sign of $(-1)^\ell = (-1)^{|\vartheta_2|}$.
This finishes the proof after dividing by the appropriate powers of $2^{-1/2}$ to rewrite the identity in terms of the rescaled Clifford multiplication $\cd$.
\end{proof}

We are now in the position to prove our the comparison theorem.

\begin{proof}[of Thm.~\ref{TheoremEqualityChPI}]
Fix $\vartheta_1, \dots, \vartheta_N \in \Omega_{\T}(X)$ throughout the proof. Define elements $\omega_a \in \Omega^1(X, \End_{\Cl_n}(\Sigma))$ and $V_a \in C^\infty(\End_{\Cl_n}(\Sigma))$ by
\begin{equation} \label{SpecialCaseOmegaV}
  \omega_a[v] \defeq - \cd(\iota_v \vartheta_a^\prime), \qquad V_a \defeq  \cd(\vartheta_a^{\prime\prime}), \qquad a =1, \dots, N.
\end{equation}
Observe that for this choice of $\omega_a$ and $V_a$, we have the equality
\begin{equation} \label{EqualityOfQs}
  Q_1 = \tilde{q}(\vartheta_1, \dots, \vartheta_N), 
\end{equation}
where $Q_\bullet$ is the stochastic process defined in \eqref{DefinitionGeneralQ} and the right hand side is the iterated Stratonovich integral \eqref{DefinitionQ} apearing in the definition of the path integral. 

On the other hand, we will now show that for this choice of $\omega_a$ and $V_a$, the operator $\Phi_t(\vartheta_1, \dots, \vartheta_N)$ from \eqref{QuantizationMap} is equal to the operator $\Psi_t^{(0)}$ appearing in Prop.~\ref{PropFeynmanKac}. 
This will follow from comparing the recursive relation \eqref{PhiTRecursion} for $\Phi_t(\vartheta_1, \dots, \vartheta_N)$ with the corresponding relation \eqref{RecursiveGeneralPhiT} for the $\Psi_t^{(a)}$.
We see that $\Phi_t = \Psi_t^{(0)}$ follows if we can show that $A_a = -F(\vartheta_a)$ and $A_{a, a+1} = - F(\vartheta_a, \vartheta_{a+1})$ for each $a$.
To this end, we calculate
\begin{equation*}
\begin{aligned}
 A_a &= \omega_a^\sharp \nabla + \frac{1}{2} \div\omega_a + V_a = - \sum_{j=1}^n \cd(\iota_{e_i} \vartheta_a^\prime) \nabla_{e_i} + \frac{1}{2}\cd(d^* \vartheta_a) +  \cd(\vartheta^{\prime\prime}_a) = -F(\vartheta_a),
\end{aligned}
\end{equation*}
where in the last equality, we used \eqref{Comparison1} to compare with the definition \eqref{DefinitionF}. 
In case that either $\vartheta_a^{\prime}$ or $\vartheta_{a+1}^{\prime}$ has degree less or equal to one, we have
\begin{equation*}
\begin{aligned}
 A_{a, a+1} &= \frac{1}{2} \omega_a^\sharp \omega_{a+1} = \frac{1}{2} \sum_{j=1}^n\cd(\iota_{e_j}\vartheta_a) \cdot \cd(\iota_{e_j} \vartheta_{a+1}) -F(\vartheta_a, \vartheta_{a+1}).
 \end{aligned}
\end{equation*}
Here we used \eqref{Comparison2} to compare with \eqref{DefinitionF} and the observation that the factor of $\frac{1}{2}$ is absorbed in the notation $\cd(\vartheta) = 2^{-|\vartheta|/2} \cc(\vartheta)$.

Using our Feynman-Kac-formula, Prop.~\ref{PropFeynmanKac} and \eqref{EqualityOfQs}, we therefore obtain that 
\begin{equation*}
  \mathbb{E}_x\bigl[S \cdot \tilde{\qq}\bigl(\rho(\vartheta_1, \dots, \vartheta_N)\bigr)\psi(\xx_1)\bigr]  = \bigl(\Phi_1(\vartheta_1, \dots, \vartheta_N) \psi\bigr)(x).
\end{equation*}
Taking conditional expectations and pointwise supertraces, followed by integration over $X$, the left hand side is just $\INT(\rho(\vartheta_1, \dots, \vartheta_N))$, while the right hand side becomes the super trace of $\Phi_1(\vartheta_1, \dots, \vartheta_N)$, represented in terms of the kernel formula \eqref{OperatorSupertrace} for the supertrace. This is precisely the Chern character $\Ch_{\Dirac}(\vartheta_1, \dots, \vartheta_N)$ as defined in \eqref{DefinitionChernCharacter}, hence the proof is finished.
\end{proof}

\section{Bismut Chern characters} \label{SectionBismutChernCharacters}

In this section, we give an application of our path integral formula to the calculation of the path integral of Bismut-Chern characters and comment on the connection to the Atiyah-Singer index theorem.

Recall definition \eqref{OmegaHat} of the complex $\widehat{\Omega}(\L X)$, respectively the $\T$-invariant subcomplex $\widehat{\Omega}(\L X)^\T$, for a Riemannian manifold $X$. The extended iterated integral map $\rho$ defined by \eqref{ExtendedIteratedIntegralMapDom} is a chain map from the entire cyclic complex $\mathsf{B}_\epsilon^\natural(\Omega_{\T}(M))$ to the complex $\widehat{\Omega}(\L X)^{\T}$ of $\T$-invariant forms. Denote its image by $\widehat{\Omega}_{\mathrm{int}}(\L X)^{\T}$.


Maybe the most prominent example of differential forms contained in the complex $\widehat{\Omega}(\L X)^{\T}$ are the Bismut-Chern-characters, first introduced in \cite{Bismut1} in the even case and \cite{Wilson} in the odd case. In fact, they can be represented by iterated integrals, hence lie in the subcomplex $\widehat{\Omega}_{\mathrm{int}}(\L X)^{\T}$.

\begin{definition}[Bismut-Chern-characters] $-$
\begin{enumerate}[(1)]
\item Let $E$ be a vector bundle with connection over $X$ and let $R$ be its curvature. The {\em Bismut-Chern-character} associated to this data is the differential form $\Ch(E) \in \widehat{\Omega}^+(\L X, \C)$ given by the formula
\begin{equation} \label{FormulaBismutChernCharacter}
  \Ch(E) = \sum_{N=0}^\infty (-1)^N \int_{\Delta_N} \tr_E\left( [\gamma\|_{\tau_N}^1]^E \bigwedge_{a=1}^N R(\tau_a) [\gamma\|_{\tau_{a-1}}^{\tau_a}]^E\right) \dd \tau.
\end{equation}
\item Let $m \in \N$ and let $g: X \rightarrow U(m)$ (the $m$-th unitary group) be a smooth map. Then the {\em odd Bismut-Chern character} $\Ch(g) \in \widehat{\Omega}^-(\L X, \C)$ is defined using the even Bismut-Chern character via the formula
\begin{equation*}
  \Ch(g) = - \int_0^1 \iota_{\partial_s}\Ch(\underline{\C}^m, \nabla^g) \dd s,
\end{equation*}
where $\underline{\C}^m$ is the trivial bundle on $X \times \R$, but with non-trivial connection  $\nabla^g = d+s g^{-1}dg$, $s$ being the coordinate of $\R$. 
\end{enumerate}
\end{definition}

\begin{remark}
Explicitly, formula \eqref{FormulaBismutChernCharacter} says that the degree $2N$-component $\Ch_N$ of $\Ch(E)$ at $\gamma \in \L X$ is given  by
\begin{equation*}
  \Ch_N[V_{2N}, \dots, V_{1}] = \frac{1}{2^{N}N!}\sum_{\sigma \in S_{2N}}\!\!\int_{\Delta_N} \!\!\!\tr_E \left([\gamma\|_{\tau_N}^1]^E \prod_{a=1}^N R\bigl(V_{\sigma_{2a}}(\tau_a), V_{\sigma_{2a-1}}(\tau_a)\bigr) [\gamma\|_{\tau_{a-1}}^{\tau_a}]^E\right) \dd \tau.
\end{equation*}
A similar formula can be given for the odd Chern character, see \cite[(6.1)]{Wilson}.
\end{remark}

\begin{proposition} \label{PropositionPropertiesCh1}
The Bismut-Chern characters are equivariantly closed. Moreover, if $j: X \rightarrow \L X$ is the inclusion, then $j^*\Ch(E) = \ch(E)$ and $j^*\Ch(g) = \ch(g)$, where $\ch(E) \in \Omega^+(X, \C)$, respectively $\ch(g) \in \Omega^-(X, \C)$ are the usual even and odd Chern character forms on $X$, defined by
\begin{equation} \label{UsualChernCharacters}
   \ch(E) = \sum_{N=0}^{n/2}\frac{(-1)^N}{N!} \tr(R^j), \qquad \ch(g) = \sum_{N=0}^{\lfloor \frac{n-1}{2}\rfloor} \frac{N!}{(2N+1)!} \tr \bigl((g^{-1} dg)^{2N+1}\bigr).
\end{equation}
\end{proposition}

This result is proven in \cite[Thm.~3.5]{Bismut1}, respectively \cite[Thm.~6.2]{Wilson}. For us, the following proposition is particularly important, which states that $\Ch(E)$ and $\Ch(g)$ are contained in the domain of the path integral map.

\begin{proposition} \label{PropositionPropertiesCh2}
Both the even and the odd Chern characters can be represented by entire, extended, iterated integrals and hence are contained in $\widehat{\Omega}_{\mathrm{int}}(\L X)^{\T}$.
\end{proposition}

\begin{proof}
In the even case, this has been shown by Getzler, Jones and Petrack \cite[\S6]{GJP}, whose construction we recall here. To write $\Ch(E)$ as an iterated integral, choose a complementary bundle $E^\perp$, meaning that $E \oplus E^\perp = \underline{\C}^m$, a trivial bundle, and let $q$ be the projection onto $E$. It is possible to arrange this in such a way that the connection on $E$ is given by the formula $\nabla = q(d + A)$, where $A = (2q-1)dq \in \Omega^1(X, \mathrm{Mat}_m(\C))$ is a matrix-valued one-form \cite{NarasimhanRamanan}. The curvature is then the two-form $R = dA + A^2$, and we set $\mathfrak{R} := A + dt \wedge R$, an element of $\Omega_\T(X, \mathrm{Mat}_m(\C))$. The Chern character is then given by $\Ch(E) = \rho(\widetilde{\Ch}(q))$, with \footnote{The formula here differs slightly from those in \cite{GJP, GueneysuLudewig} owing to the use of the bar complex instead of the cyclic complex of Connes used there.}
\begin{equation} \label{BismutChernIteratedIntegral}
  \widetilde{\Ch}(q) = \sum_{N=0}^\infty (-1)^N \sum_{m=0}^N\tr \bigl(\underbrace{\mathfrak{R}, \dots, \mathfrak{R}}_m, dt \wedge q, \underbrace{\mathfrak{R}, \dots, \mathfrak{R}}_{N-m}\bigr),
\end{equation}
where the trace is taken as in \cite[Def.~7.2]{GueneysuLudewig}; see Thm.~6.5 in \cite{GJP} for a proof of this claim. It is easy to see that $\widetilde{\Ch}(q)$ satisfies the necessary estimates in order to have finite entire norms \eqref{EntireSeminorms}; compare also \cite[\S 7]{GueneysuLudewig}. 

A similar formula can be given for the odd Bismut-Chern characters. Explicitly, given a smooth map $g: X \rightarrow \mathrm{U}_k$ and $s \in \R$, set
\begin{equation}
\begin{aligned} \label{AandRodd}
A_s &= sg^{-1}dg \\
R_s &= dA_s + A_s^2 = -s(1-s) g^{-1} dg \wedge g^{-1} dg.
\end{aligned}
\end{equation}
and $\mathfrak{R}_s = A_s + dt \wedge R_s \in \Omega_\T(X)$. Then the formula for the ``combinatorial'' odd Chern character is
\begin{equation} \label{BismutChernIteratedIntegralOdd}
  \widetilde{\Ch}(g) = -\sum_{N=0}^\infty (-1)^N \sum_{m=0}^N \int_0^1 \tr (\underbrace{\mathfrak{R}_s, \dots, \mathfrak{R}_s}_m, dt \wedge g^{-1}dg, \underbrace{\mathfrak{R}_s, \dots, \mathfrak{R}_s}_{N-m}) \dd s;
\end{equation}
 this has been shown in \cite[\S5]{CacciatoriGueneysu}.
\end{proof}

\begin{remark} \label{InfiniteBismutChern}
The Chern character $\Ch(q)$ can be written elegantly using the shuffle product. In fact, since  in the grading-shifted space $\Omega_{\T}(X)[1]$, $A$ has degree zero and $R$ has degree $2$, the degree $2\ell$ part of $\Ch(q)$ is the infinite sum of shuffle products
\begin{equation*}
  \sum_{N=0}^\infty (dt \wedge p) \shuffle (\underbrace{dt \wedge R, \dots, dt \wedge R}_{\ell}) \shuffle (\underbrace{A, \dots, A}_N).
\end{equation*}
Hence even in each fixed degree, $\Ch(q)$ is an infinite sum and the entire topology \eqref{EntireSeminorms} is needed to make sense of it. In particular, the Continuity Lemma~\ref{LemmaIteratedIntegralMapContinuous} is needed to apply the iterated integral map to $\tilde{\Ch}(q)$, a fact that is omitted in the discussion of \cite[\S 6]{GJP}.
\end{remark}

\paragraph{Application to index theory.} 
To close this section, we calculate the path integrals of the Bismut-Chern characters, which are given in terms of an index, respectively a spectral flow.

\begin{proposition} \label{PropFormulasBC}
Let $X$ be a compact spin manifold. Then the path integrals of the Bismut-Chern characters are given as follows.
\begin{enumerate}[{\normalfont(a)}]
\item If $\Dirac_E$ is the Dirac operator twisted by $E$, then
\begin{equation} \label{IndexTheoremEven}
\INT\bigl(\Ch(E)\bigr) = i^{n/2}\mathrm{ind}(\Dirac_E).
\end{equation}
\item Let $\Dirac$ be the Dirac operator acting on $\Sigma \otimes \underline{\C}^k$ and denote by  $\mathrm{sf}(\Dirac, g^{-1} \Dirac g)$ the spectral flow between the Dirac operators $\Dirac$ and $g^{-1} \Dirac g$, where $g : C \to U(m)$. Then 
\begin{equation} \label{IndexTheoremOdd}
\INT\bigl(\Ch(g)\bigr) = -(-i)^{\frac{n+1}{2}} \sqrt{2\pi} \cdot \mathrm{sf}(\Dirac, g^{-1} \Dirac g).
\end{equation}
\end{enumerate}
\end{proposition}

It follows directly from the definition of the supertrace that $\str$ is and even respectively odd, in even and respectively odd dimensions. This implies that $\INT$ has the same property. As the spectral flows of Dirac operators vanish in even dimensions (by symmetry of the spectrum), \eqref{IndexTheoremOdd} is trivial in even dimensions, while \eqref{IndexTheoremEven} is trivial in odd dimensions (there, all indices vanish).

\begin{proof}
Let $E$ be a Hermitean vector bundle with compatible connection $\nabla$ over $X$. As in the proof of Proposition \ref{PropositionPropertiesCh2}, realize $E = \mathrm{im}(q)$ for some projection $q \in C^\infty(X, \mathrm{Mat}_m(\C))$. Let $A$ be the corresponding connection 1-form and $R = dA + A^2$ its curvature. 

As in \S\ref{SectionComparison}, let $\UU_\bullet$ be the horizontal lift of the Brownian motion $\xx_\bullet$ to the spin principal bundle $P$, and let $\mathbf{A} \in \Omega^1(P, \mathrm{Mat}_m(\C))^{\Spin_n}$ and $\mathbf{R} \in \Omega^2(P, \mathrm{Mat}_m(\C))^{\Spin_n}$ be the horizontal lifts of $A$, respectively $R$. Let $\bfQ_t$ be the solution to the $\C_n \otimes \mathrm{Mat}_m(\C)$-valued Stratonovich differential equation
\begin{equation*}
  \dd \bfQ_t = \bfQ_t \cdot\bigl( \mathbf{A} (\UU_t) * \dd \UU_t - \cd\bigl(\mathbf{R}(\UU_t)) \dd t  \bigr), \qquad \bfQ_0 = \mathrm{id}.
\end{equation*}
Here, $\cd(\mathbf{R}) = \frac{1}{2}\sum_{i < j} \cc(e_i)\cc(e_j) \otimes \mathbf{R}[e_i, e_j] \in C^\infty(P, \Cl_n \otimes \mathrm{Mat}_m(\C))^{\Spin_n}$ is the term appearing in the general Lichnerowicz formula for the twisted Dirac operator (\cite[Thm.~3.52]{BGV}).
Then $\bfQ_t$ can be written as the iterated Stratonovich integral 
\begin{equation*}
\bfQ_t = \sum_{N=0}^\infty (-1)^N \int_0^t\int_0^{\tau_N} \cdots \int_0^{\tau_2} \prod_{a=1}^N \bigl(\cd\bigl(\mathbf{R}(\UU_t)\bigr) \dd t -  \mathbf{A} (\UU_t) * \dd \UU_t   \bigr).
\end{equation*}
In view of \eqref{DefinitionQ} and the iterated integral expression \eqref{BismutChernIteratedIntegral} of $\Ch(E)$, we therefore obtain 
\begin{equation*}
 \tilde{q}\bigl(\Ch(E))=  \int_0^1 \UU_1 \tr\bigl(\bfQ_t \, q \,\bfQ_{1-t}\bigr) \UU_0^* \dd t.
\end{equation*}
By the vector-valued Feynman-Kac formula (see \cite{MR2735286}) and the Lichnerowicz formula $\Dirac_A^2 = \nabla^*\nabla + \frac{1}{4} \mathrm{scal} + \cc(R)$ for the Dirac operator $\Dirac_A$ twisted with the connection $\nabla^\Sigma \otimes (d+A)$ on $\Sigma \otimes \underline{\C}^m$, we have by cyclicity of the supertrace
\begin{equation*}
\begin{aligned}
  \INT(\Ch(E)) = \int_0^1 \Str \otimes \mathrm{Tr}_{\underline{\C}^m}\bigl(e^{-t\Dirac_A/2} q e^{-(1-t)\Dirac_A/2}\bigr) \dd t &= \Str  \otimes \mathrm{Tr}_{\underline{\C}^m}(q e^{-\Dirac_A/2}) \\
  &= i^{n/2}\Str_{\Sigma_\C \otimes E}(e^{-\Dirac_E/2}),
\end{aligned}
\end{equation*}
where we observed $q e^{-\Dirac_A/2} = e^{-\Dirac_E/2}$ and used \eqref{TraceComparisonEven} to convert the real supertrace to the complex one. Now for the complex supertrace, the McKean-Singer formula applies (see e.g.\ \cite[Thm.~3.50]{BGV}), giving statement (a).

To see \eqref{IndexTheoremOdd}, we may assume that $n = \dim(X)$ is odd. 
Here, similar to the above, one obtains that $\INT(\Ch(g))$ equals
\begin{equation*}
- \int_0^1 \Str_{\Sigma \otimes \underline{\C}^m} \Bigl(\cd(g^{-1} dg) e^{-\Dirac_s^2/2}\Bigr) \dd s = - (-i)^{\frac{n+1}{2}}\int_0^1 \mathrm{Tr}_{\Sigma_\C \otimes \underline{\C}^m} \Bigl(\cc(g^{-1} dg) e^{-\Dirac_{s}^2/2}\Bigr) \dd s.
\end{equation*}
In the above equality, we used that since $n$ is odd, for an \emph{odd} element $a \in \Cl_n$, we have $2^{-1/2} \str(a) = (-i)^{\frac{n+1}{2}}\tr_\C(a)$, by \eqref{TraceComparisonOdd} and \eqref{DefinitionSupertrace}. 
This swallows the additional factor of $2^{-1/2}$ coming from the use of the rescaled Clifford multiplication there.
From Getzler \cite[Corollary~2.7]{GetzlerOdd}, we have
\begin{equation} \label{GetzlerIndexTheorem}
  \mathrm{sf}(\Dirac, g^{-1} \Dirac g) = (2\pi)^{-1/2} \int_0^1 \mathrm{Tr}_{\Sigma_\C \otimes \underline{\C}^k} \Bigl( \dot{\Dirac}_{s} e^{-\Dirac_{s}^2/2} \Bigr) \dd s.
\end{equation}
The result now follows after noticing that $\dot{\Dirac}_s = \cc(g^{-1} dg)$.
\end{proof}

\begin{remark}
Observe that since both $\Ch(E)$ and $\Ch(g)$ are such that Thm.~\ref{TheoremEqualityChPI} applies.
Hence we have
\begin{equation*}
  \INT\bigl(\Ch(E)\bigr) = \Ch_{\Dirac}\bigl(\tilde{\Ch}(q)\bigr), 
\end{equation*}
i.e., we can also compute the path integral of the Bismut-Chern characters using the algebraic Chern character. In fact, this gives another proof of Prop.~\ref{PropFormulasBC}(a), using the results of \cite[\S8]{GueneysuLudewig}.
\end{remark}

Using Thm.~\ref{TheoremEqualityChPI}, the localization formula \eqref{LocalizationFormula} for $\Ch_{\Dirac}$ can be carried over to a localization formula for $\INT$, which reads as follows and is the loop space analog of the finite-dimensional localization principle of Duistermaat-Heckmann \cite{DuistermaatHeckmann, BerlineVergne}. 

\begin{theorem} \label{ThmLocalizationI}
For any iterated integral $\theta \in \Omega_{\mathrm{int}}(\L X)$ that is equivariantly closed and satisfies the assumptions of {\normalfont Thm.~\ref{TheoremEqualityChPI}}, we have 
\begin{equation} \label{EqLocalizationI}
  \INT(\theta) = (2 \pi)^{-n/2} \int_X \hat{A}(X) \wedge j^*\theta,
\end{equation}
where $j: X \rightarrow \L X$ is the inclusion map. 
\end{theorem}

Observe that $\Ch(E)$ satisfies the assumptions of Thm.~\ref{ThmLocalizationI}.
Since $\Ch(E)$ is equivariantly closed, we can alternatively use \eqref{EqLocalizationI} in order to calculate the iterated integral of the Bismut-Chern characters. 
In view of \eqref{UsualChernCharacters}, in the even case, the result  is
\begin{equation*}
  \INT(\Ch(E)) = (2 \pi)^{-n/2}\int_X \hat{A}(X) \wedge \ch(E).
\end{equation*}
Together with Prop.~\ref{PropFormulasBC}, this is the Atiyah-Singer index formula. 
A similar result can be obtained in the odd case.

\bibliography{LiteraturLoopSpaces}

\end{document}